\documentclass[preprint,12pt]{article}
\usepackage{amssymb}
\usepackage{mathrsfs}
\usepackage{amsfonts}
\usepackage{graphicx}
\usepackage{colortbl,dcolumn}
\usepackage{amsmath}
\usepackage{psfrag}
\usepackage{booktabs}
\usepackage{cite}
\numberwithin{equation}{section}
\usepackage[top=1.1in, bottom=1.2in, left=0.8in, right=0.8in]{geometry}

\newcommand{\R}{\mathbb{R}}
\newcommand{\N}{\mathbb{N}}
\newcommand{\E}{\mathbb{E}}

\newcommand{\dd}{\text{d}}

\newtheorem{thm}{Theorem}[section]

\newtheorem{lem}[thm]{Lemma}
\newtheorem{prop}[thm]{Proposition}
\newtheorem{cor}[thm]{Corollary}
\newtheorem{rem}[thm]{Remark}

\newtheorem{assumption}[thm]{Assumption}

\begin{document}
\title{Sharp mean-square regularity results for SPDEs with fractional noise and optimal convergence rates for the numerical approximations \footnotemark[1]}

\author{
Xiaojie Wang$\,^\text{a}$ \footnotemark[2] 
\quad 
Ruisheng Qi$\,^\text{b}$   
\quad 
Fengze Jiang$\,^\text{c}$ 
\\
\footnotesize $\,^\text{a}$ School of Mathematics and Statistics, Central South University, Changsha, China\\
\footnotesize $\,^\text{b}$ School of Mathematics and Statistics, Northeastern University at Qinhuangdao, Qinhuangdao, China\\
\footnotesize $\,^\text{c}$ School of Mathematics and Statistics, Huazhong University of Science and Technology, Wuhan, China\\
}

       \maketitle

       \footnotetext{\footnotemark[1] X.W. was partially supported by NNSF of China (No.11301550, No.11171352).
                R.Q. was supported by  Research Fund for Northeastern University at Qinhuangdao (No.XNB201429),
                Fundamental Research Funds for Central Universities  (No.N130323015),  Science and Technology 
                Research Founds for Colleges and Universities in Hebei Province (No.Z2014040), 
                Natural Science Foundation of Hebei Province, China (No.A2015501102).
                }
        \footnotetext{\footnotemark[2] x.j.wang7@gmail.com, x.j.wang7@csu.edu.cn}
       \begin{abstract}
          {\rm\small 
           This article offers sharp spatial and temporal mean-square regularity results for a class of 
           semi-linear parabolic stochastic partial differential equations (SPDEs) driven by infinite 
           dimensional fractional Brownian motion with the Hurst parameter greater than one-half.
           In addition, mean-square numerical approximation of such problem are investigated, 
           performed by the spectral Galerkin method in space and the linear implicit Euler method in time.
           The obtained sharp regularity properties of the problems enable us to identify optimal mean-square
           convergence rates of the full discrete scheme. These theoretical findings are accompanied by 
           several numerical examples.
           } \\


\end{abstract}

\section{Introduction}\label{sec:introduction}

Numerical analysis of evolutionary stochastic partial differential equations (SPDEs) is currently an active area 
of research and there has been an extensive literature on numerical methods for SPDEs driven by
standard (possibly cylindrical) $Q$-Wiener process
\cite{kruse2014strong,kruse2012optimal,kruse2014optimal,lord2014introduction,jentzen2009numerical,jentzen2011taylor}, to just mention a list of comprehensive references. 
Although the theoretical analysis of SPDEs driven by infinite dimensional fractional Brownian motion (fBm in short)
has attracted increasing attention in the past decades (see, e.g., \cite{duncan2002fractional,duncan2009semilinear,grecksch1999parabolic,maslowski2003evolution,pasik2006linear,hu2009stochastic,sanz2009mild} and references therein), 
they are still not well-understood, especially from the numerical point of view.
The purpose of the present work is to provide sharp mean-square regularity properties of SPDEs with 
infinite dimensional fractional noise and to identify optimal mean-square convergence rates 
of the numerical approaximations.


Given a real separable Hilbert space $(V, \langle \cdot, \cdot \rangle, \|\cdot\| )$ with 
scalar product $\langle \cdot, \cdot \rangle$ and norm $\|\cdot\| = \langle \cdot, \cdot \rangle^{\frac{1}{2}}$,
we let $A \colon \mathcal{D}(A) \subset V \rightarrow V$ be  a densely defined, linear unbounded, 
positive self-adjoint operator with compact inverse.
Let $\left(\Omega,\mathcal {F},\mathbb{P}\right)$ be a probability space 
and let $\{W^H (t) \}_{t \in [0, T] }$ be a standard cylindrical fBm 
with Hurst parameter $H \in (\tfrac12, 1)$, defined by the following formal series
\begin{equation} \label{eq:intro-Wiener-representation}
W^H (t) := \sum_{n = 1}^{\infty} w_n^H ( t ) e_n,
\quad
t \in [0, T],
\end{equation}
where $\{w_n^H ( t )\}_{n \in \N}, t \in [0, T]$ are a sequence of independent real-valued standard fBm 
each with the same Hurst parameter $H \in ( \tfrac12, 1)$ and $\{e_n\}_{n \in \N}$ being a complete orthonormal basis of $V$.
Now let us consider the following semilinear parabolic SPDEs
in $V$, driven by an infinite dimensional fractional Brownian motion (fBm),
\begin{equation}\label{eq:spde}
\begin{split}
\left\{
    \begin{array}{lll} \dd X(t) + A X(t)\, \dd t = F ( X(t) ) \,\dd t +  \Phi \, \dd W^H(t), \quad  t \in (0, T], \\
     X(0) = \xi,
    \end{array}\right.
\end{split}
\end{equation}
where $F \colon V \rightarrow V$ and $\Phi \colon V \rightarrow V$ are deterministic mappings.
Under certain assumptions specified later,  particularly including
\begin{equation}\label{eq:intro-AQ-condition}
\|A^{\frac{\beta-1}{2}} \Phi \|_{\mathcal{L}_2(V) }  < \infty, \quad \text{ for some } \beta \in ( 1 - 2 H, 1],
\end{equation}
\eqref{eq:spde} admits a unique mild solution 
$X: [0,T] \times \Omega \rightarrow V$ with continuous sample path, given by
\begin{equation}\label{eq:intro-mild}
    X(t)
    =
    E(t) \xi
    +
    \int_0^t E(t-s) F ( X( s ) ) \, \dd s
    +
    \int_0^t E(t-s) \Phi \, \dd W^H(s), \quad  \mathbb{P} \mbox{-a.s.}.
\end{equation}
Here $E(t) = \exp(-t A), t \geq 0$ represents an analytic semigroup generated by $-A$. 
As shown in the main regularity result, Theorem \ref{thm:SPDE-regularity-result}, 
the mild solution possesses the following regularity properties:
\begin{equation}\label{eq:intro-space-regularity}
\sup_{t \in [0, T]} \|  X(t) \|_{L^2 ( \Omega; V_{2H + \beta - 1}) } < \infty,
\end{equation}
and
\begin{equation}\label{eq:intro-time-regularity}
\sup_{t \neq s} \frac{ \|  X(t) - X(s) \|_{L^2 ( \Omega; V_{\delta}) } } 
{(t -s)^{ \frac{ 2H + \beta - 1 - \delta }{2} } }
< \infty,
\quad
\delta \in [0, 2H + \beta - 1 ].
\end{equation}
Clearly, the parameter $\beta$ in \eqref{eq:intro-AQ-condition} used to characterize the spatial regularity 
of the operator $\Phi$ also determines the spatial and temporal regularity of the mild solution \eqref{eq:intro-mild}.
If $\Phi \in \mathcal{L}_2(V)$, which corresponds to the trace-class noise with the fBm in \eqref{eq:spde}
being of covariance type, then \eqref{eq:intro-AQ-condition} is satisfied with $\beta = 1$ and the mild solution
$X(t), t \in [0, T] $ takes values in $L^2 ( \Omega; V_{2H})$ and is mean-square H\"older continuous in $V_{\delta}$
with H\"older exponent $H - \frac{ \delta }{2}$ for $\delta \in [0, 2H ]$.
For an interesting case when it is only assumed that $\Phi \in \mathcal{L}(V)$ ($\Phi = I_V$ for example), 
\eqref{eq:intro-AQ-condition}  is fulfilled with $\beta < \tfrac12$ in  one space dimension, 
$\beta < 0$ in two space dimensions and  $\beta < -\tfrac12$ in three dimensions 
\cite[section 6.1]{kovacs2008introduction}.  Therefore, regularity results 
\eqref{eq:intro-space-regularity}-\eqref{eq:intro-time-regularity} tell us that, 
distinct from parabolic SPDEs driven by standard cylindrical $I$-Wiener process 
($H = \tfrac12$), SPDEs driven by standard cylindrical fBm ($ \Phi = I_V$) allow for a mild solution 
with a positive order of regularity in multiple spatial dimensions.
Also, it is worthwhile to point out that, the sharp regularity results for the mild solution \eqref{eq:intro-mild} 
is credited to sharp regularity results of the stochastic convolution. 
A key tool utilized to analyze the sharp regularity properties of the stochastic convolution is 
Lemma \ref{lem:optimal-regularity-lemma}, the proof of which is due to a very careful use of the smoothing property of the analytic semigroup. 

Additionally, in this paper we study a full discretization of \eqref{eq:intro-mild} via a spectral 
Galerkin method for the spatial approximation and the linear implicit Euler method 
for the temporal discretization. By $X(t_m)$ we denote the mild solution 
\eqref{eq:intro-mild} taking values at the temporal grid points 
$t_m = m \tau$ and by $\bar{ X}^N_{m}$ the numerical approximations of $X(t_m)$, 
produced by the proposed fully discrete scheme with the time step-size $\tau > 0$. 
Based on the above sharp regularity properties, one can measure the resulting 
approximation error as follows (Corollary \ref{cor:full-discrete-scheme-error-bound}):
\begin{equation} \label{eq:Intro-numerical-main-result}
\| X ( t_{m} ) -  \bar{ X}^N_{m} \|_{L^2 ( \Omega; V) }
     \leq
C ( 1 + \| \xi \|_{L^2 ( \Omega; V_{2H + \beta - 1}) })
\big( \lambda_{N + 1}^{- \frac{ 2 H + \beta - 1}{2} } + \tau^{ \frac{ 2 H + \beta - 1 }{2} } \big),
\end{equation}
where $\{\lambda_n\}_{n \in \N }$ are the eigenvalues 
of the linear operator $A$. Comparing this with the regularity results \eqref{eq:intro-space-regularity} 
and \eqref{eq:intro-time-regularity}, one can easily observe, 
the obtained convergence rates in \eqref{eq:Intro-numerical-main-result} are optimal in the sense that 
the orders of convergence in space and in time agree with the order of the spatial and temporal regularity 
of the mild solution, respectively. It must be emphasized that the derivation of 
\eqref{eq:Intro-numerical-main-result} is not an easy task and requires
a variety of delicate error estimates, which are elaborated in section \ref{sec:strong.convergence}.
Unlike the case of standard cylindrical Wiener process, It\^o's isometry for the stochastic integral 
with respect to a cylindrical fBm gets involved with a multiple integral (cf. \eqref{eq:intro-ito-isometry} below).
This essentially makes the treatments of the corresponding error terms, i.e., $|J_3|^2$ in \eqref{eq:J3} and 
$|\mathbb{J}_3|^2$ in \eqref{eq:J3-estimate},  significantly more difficult and demanding.

Before closing this introduction section, we recall a few closely relevant works. 
In \cite{duncan2002fractional,duncan2009semilinear}, solutions of linear and semi-linear SPDEs 
with additive fractional noises are investigated, but sharp mean-square regularity results were
not reported there. Indeed, we must admit that,  our work is inspired by \cite{kruse2012optimal} 
and \cite{kruse2014optimal},  where the former established sharp mean-square 
regularity properties of parabolic SPDEs driven by standard $Q$-Wiener process 
of covariance type ($\beta =1$ in our setting) and the latter recovered 
optimal convergence rates of the numerical approximations.  Nevertheless, as already discussed earlier,
a more complicated form of It\^o's isometry in the fBm setting makes 
both the regularity analysis (see the proof of Lemma \ref{lem:optimal-regularity-lemma}) 
and especially the approximation error analysis (see estimates of $|J_3|^2$ and $|\mathbb{J}_3|^2$ 
in section \ref{sec:strong.convergence}) much more involved and new techniques are needed.
At last, we would like to mention two recent publications \cite{cao2015well,cao2016approximating}, 
where Cao, Hong and Liu examined strong approximations of various SPDEs driven by  fractional noise 
with $H \in (0, \tfrac12)$. Instead of the semigroup approach adopted in this paper,  
they used the Green function framework.

The rest of this paper is organized as follows. In the next section we collect some basic facts and
define the stochastic integral with respect to the infinite dimensional fBm. 
Section \ref{sec:optimal-regularity} is devoted to sharp regularity analysis 
of the underlying SPDEs under standard assumptions.  With the sharp regularity results, 
we show optimal mean-square convergence rates of the numerical approximations
in section \ref{sec:strong.convergence}. 
Numerical results are included in section \ref{sec:numerical-result} 
to test previous theoretical findings.
Finally, the paper is concluded with further comments.

\section{Preliminaries} \label{sec:prelim}

On a real separable Hilbert space $(V, \langle \cdot, \cdot \rangle, \|\cdot\| )$, 
by $\mathcal{L}(V)$ we denote the space of bounded linear operators from 
$V$ to $V$ endowed with the usual operator norm $\| \cdot \|_{\mathcal{L}(V)}$.
Additionally, we denote by $\mathcal{L}_2(V) \subset \mathcal{L}(V)$ the subspace
consisting of all Hilbert-Schmidt operators from $V$ to $V$ \cite{da2014stochastic}. 
It is known that $\mathcal{L}_2(V)$ 
is a separable Hilbert space, equipped with the scalar product and norm
\begin{equation}\label{eq:defn-L2}
\langle \Gamma_1, \Gamma_2  \rangle_{\mathcal{L}_2(V)} 
: = 
\sum_{n \in \N} \langle \Gamma_1 \eta_n, \Gamma_2 \eta_n  \rangle,
\quad
\| \Gamma \|_{ \mathcal{L}_2(V) } 
:=
\Big (
\sum_{n \in \N} \| \Gamma \eta_n \|^2
\Big )^{\frac12},
\end{equation}
independent of the particular choice of ON-basis $\{\eta_n\}_{n \in \N}$ of $V$.
Below we sometimes write $\mathcal{L}_2 := \mathcal{L}_2(V)$ for brevity.
If $\Gamma \in \mathcal{L}(V)$ and $\Gamma_1, \Gamma_2 \in \mathcal{L}_2(U)$, then
the following inequalities hold,
\begin{equation}\label{eq:L-L2-norm-inequality}
| \langle \Gamma_1, \Gamma_2  \rangle_{\mathcal{L}_2(V)} | 
\leq 
\| \Gamma_1 \|_{\mathcal{L}_2(V)} 
\| \Gamma_2  \|_{\mathcal{L}_2(V)},
\quad
\| \Gamma \Gamma_1\|_{\mathcal{L}_2(V)} \leq \|\Gamma\|_{\mathcal{L}(V)} \|\Gamma_1\|_{\mathcal{L}_2(V)}.
\end{equation}
As a main target of this section, we are to define 
a stochastic integral of the form
\begin{equation}\label{eq:Stoch-Integeral-defn}
I ( \Psi; T ) := \int_0^T \Psi (s) \, \text{d} W^H (s), \quad H \in (\tfrac12, 1),
\end{equation}
for an integrand $\Psi \colon [0, T] \rightarrow \mathcal{L}_2(V) $, where $W^H$ is the cylindrical fBm represented by \eqref{eq:intro-Wiener-representation}.
There are a variety of ways to define such stochastic integral in existing literature. 
Here we adopt a transparent approach used in \cite{duncan2002fractional}
and recall an inequality of integral form as \cite[Lemma 2.1]{duncan2002fractional}.
\begin{lem}\label{lem:double-intefral-inequality}
Let $\chi \in L^p (0, T; \R)$ for $p > \tfrac{1}{H}$ be a deterministic function. 
Then there exists a constant $C_T > 0$ only depending on $T$ such that
\begin{equation}
\int_0^T \! \int_0^T \chi (u) \chi (v) \phi (u - v) \dd u \dd v \leq C_{T} \| \chi \|^2_{L^p (0, T; \R)},
\end{equation}
where and below for simplicity of presentation we denote
\begin{equation}
\phi( y ) :=  \alpha_H | y |^{2 H - 2}, \ y \in \R  \quad \text{ with } \quad  \alpha_H = H ( 2 H - 1). 
\end{equation}
\end{lem}
Then for a deterministic $V$-valued function $g \in L^p (0, T; V)$ with $p > 1/ H$ 
and a scalar fractional Brownian motion $w^H$,  we first define the stochastic integral
\begin{equation} \label{eq:stoch-integral-scalar}
\int_0^T g (s)\, \text{d} w^H(s).
\end{equation} 
The scalar fractional Brownian motion $w^H$ for $H \neq \tfrac12$ is 
neither a Markov process nor a semi-martingale, 
but a centered Guassian process with continuous samples and covariance function
\begin{equation}
R_H ( s , t ) = \tfrac12 ( s^{2H} + t^{2H} - | t - s |^{ 2H } ).
\end{equation}
More details of the scalar fBm can be found in \cite{biagini2008stochastic}.
At the moment let us focus on the definition of the stochastic integral \eqref{eq:stoch-integral-scalar}.
Let $\Xi$ be the family of $V$-valued step functions, defined by
\begin{equation}
\Xi := \Big \{ 
g \colon g (s) = \sum_{i = 0}^{n - 1} g_i 1_{ [ t_i, t_{i+1} ) } (s), \, 0 = t_0 < t_1 < \cdots < t_n = T, \,  g_i \in V 
\Big \}.
\end{equation}
For such special function $ g \in \Xi$ we define the stochastic integral \eqref{eq:stoch-integral-scalar} as
\begin{equation}
\int_0^T g (s)\, \text{d} w^H(s) : =   \sum_{i = 0}^{n - 1} g_i \big(  w^H (t_{i+1}) - w^H (t_{i}) \big).
\end{equation}
It is not difficult to check that, the mean of this random variable is zero and its second moment satisfies
\begin{equation} \label{eq:scalar-fBm-Ito-isometry}
\E \Big[ \Big \| \int_0^T g (s)\, \text{d} w^H(s) \Big \|^2 \Big] = \int_0^T \int_0^T \langle g(u), g(v) \rangle \phi (u - v) \text{d} u   \text{d} v \leq C_{p, T} \| g \|^2_{L^p (0, T; V) }
\end{equation}
for some constant $C_{p, T}$ only depending on $p, T$. The inequality in \eqref{eq:scalar-fBm-Ito-isometry}
holds true thanks to Cauchy–Schwarz inequality and Lemma \ref{lem:double-intefral-inequality}.
Since $\Xi$ is dense in $L^p (0, T; V)$, the stochastic integral can be (almost surely)
uniquely extended from $\Xi$ to $L^p (0, T; V)$.

Now we return to the definition of  $I ( \Psi; T )$ as indicated in \eqref{eq:Stoch-Integeral-defn}. 
Assume
\begin{equation}
\Psi ( \cdot ) x \in  L^p (0, T; V)
\:\:\:
\forall \,
x \in V,
\: \text{ and }
\:
\int_0^T \int_0^T \| \Psi (u) \|_{\mathcal{L}_2 (V) } \| \Psi (v) \|_{\mathcal{L}_2 (V) } \phi ( u - v ) \text{d} u \text{d} v < \infty.
\end{equation}
Under theses assumptions, we define the stochastic integral $I( \Psi; T )$ as
\begin{equation} \label{eq:defn-stochastic-integeral}
I( \Psi; T ) = \int_0^T \Psi (s) \, \text{d} W^H (s) := \sum_{n = 1}^{\infty} 
\int_0^T  \Psi (s) e_n \, \text{d} w_n^H (s),
\end{equation}
where the summation is defined in mean square.
Since $\Psi (\cdot) e_n \in L^p (0, T; V)$ for each $n \in \N$, 
all summands in \eqref{eq:defn-stochastic-integeral} are well-defined
due to the definition of \eqref{eq:stoch-integral-scalar}
and they are mutually independent Gaussian random variables.
Easy calculations show that, the series in \eqref{eq:defn-stochastic-integeral} is a zero mean, 
$V$-valued Gaussian random variable  satisfying the following It\^o's isometry:
\begin{equation} \label{eq:intro-ito-isometry}
\begin{split}
\E \Big[ \Big \| \int_0^T \Psi (s)\, \text{d} W^H(s) \Big \|^2 \Big]  & = \sum_{n = 1}^{\infty} \E \Big[ \Big \| \int_0^T \Psi (s) e_n \, \text{d} w^H(s) \Big \|^2 \Big] 
\\ & =
\sum_{n = 1}^{\infty} \int_0^T \int_0^T \big\langle \Psi(u) e_n, \Psi (v) e_n \big\rangle \phi (u - v) \text{d} u \, \text{d} v
\\ & =
\int_0^T \int_0^T \big\langle \Psi(u) , \Psi (v) \big\rangle_{\mathcal{L}_2}  \phi (u - v) \text{d} u \,  \text{d} v,
\end{split}
\end{equation} 
which serves as an important tool in the error analysis later and where by assumption
\begin{equation}
\int_0^T \int_0^T \big\langle \Psi(u) , \Psi (v) \big\rangle_{\mathcal{L}_2}  \phi (u - v) \text{d} u   \text{d} v 
\leq
\int_0^T \int_0^T \| \Psi(u) \|_{\mathcal{L}_2 }  \|\Psi (v) \|_{\mathcal{L}_2}  \phi (u - v) \text{d} u   \text{d} v 
< \infty.
\end{equation}

\section{Sharp regularity results}
\label{sec:optimal-regularity}
This section aims to analyze mean-square regularity properties of \eqref{eq:intro-mild} in both space and time.
To begin with, we make the following assumptions
\begin{assumption}[Linear operator A]\label{ass:A}
Let $(V,\, \langle \cdot, \cdot \rangle, \,\|\cdot \|)$ be a real separable Hilbert space and let $A \colon \mathcal{D}(A) \subset V \rightarrow V$ 
be a linear, densely-defined,  positive self-adjoint unbounded operator with compact inverse.
\end{assumption}
This assumption guarantees  that $-A$ generates an analytic semigroup $E(t) = e^{-t A}, t \geq 0$ on $V$.
Furthermore, there exists an increasing sequence of real numbers $\{\lambda_i\}_{i =1}^{\infty}$ and 
an orthonormal basis $\{e_i\}_{i =1}^{\infty}$ such that $A e_i = \lambda_i e_i$ and
\begin{equation}
0 < \lambda_1 \leq \lambda_2 \leq \cdots \leq \lambda_n (\rightarrow \infty).
\end{equation}
This allows us to define the fractional powers of $A$, i.e., $A^\gamma, \gamma \in \mathbb{R}$ 
and the Hilbert space $V_{\gamma} :=\text{dom}(A^{\frac{\gamma}{2}})$, equipped with inner product  
$\langle \cdot, \cdot \rangle_\gamma : = \langle A^{\frac{\gamma}{2}} \cdot, 
A^{\frac{\gamma}{2}} \cdot \rangle$ 
and norm $\| \cdot \|_{\gamma} = \langle \cdot, \cdot \rangle_\gamma^{\frac12} $ 
\cite[Appendix B.2]{kruse2014strong}. Moreover, $V_0 = V$ and 
$V_{\gamma} \subset V_{\delta},  \gamma \geq \delta$. It is also well-known that \cite{pazy1983semigroups}
\begin{equation}
\begin{split}
\label{eq:E.Inequality}
\| A^\gamma E(t)\|_{\mathcal{L}(V)} \leq& C t^{-\gamma}, \quad t >0, \gamma \geq 0,
\\
\| A^{-\rho} (I-E(t))\|_{\mathcal{L}(V)} \leq& C t^{\rho}, \quad t >0, \rho \in [0,1],
\end{split}
\end{equation}
which together imply
\begin{equation}\label{eq:Semigroup-differ}
\| E (t) - E (s) \|_{\mathcal{L}(V)} =  \big \| E (s) \big(  E (t - s) - I \big) \big \|_{\mathcal{L}(V)} \leq C s^{-\rho} (t - s)^{\rho} , 
\quad t> s> 0, \, \rho \in [0, 1].
\end{equation}
Throughout this paper, 
by $C$ and $c_H$ we mean various constants, not necessarily the same at each occurrence, that are independent of the discretization parameters.

\begin{assumption}[Nonlinearity]\label{ass:F}
Let $F \colon V \rightarrow V $ be a deterministic mapping satisfying
\begin{align}
\label{eq:F.conditions1}
\| F( u ) \|  & \leq L (\| u \| +1), \quad\,   u \in V,
\\
\label{eq:F.conditions2}
\| F(u) - F (v) \|  & \leq L \| u - v \|, \qquad  u, v \in V
\end{align}
for some constant $L \in (0, \infty)$.
\end{assumption}

\begin{assumption}[Noise term]\label{ass:Phi}
Let $\{ W^H (t) \}_{t \in[0, T]}$ be a cylindrical fractional Brownian motion on a 
probability space  $\left(\Omega,\mathcal {F},\mathbb{P} \right)$ 
with a normal filtration $\{\mathcal{F}_t\}_{ t\in [0, T] }$, expressed by
\begin{equation} \label{eq:Wiener-representation}
W^H (t) := \sum_{n = 1}^{\infty} w_n^H ( t ) e_n,
\end{equation}
where $\{w_n^H ( t )\}_{n \in \N}$ is a sequence of independent real-valued standard fractional Brownian motions 
each with the same Hurst parameter $H \in ( \tfrac12, 1)$ and $\{e_n\}_{n \in \N}$ is a complete orthonormal basis 
of $V$.
Assume further that, the deterministic mapping $\Phi \colon V \rightarrow V$ satisfies
\begin{equation}
\label{AQ_condition}
\|A^{\frac{\beta-1}{2}} \Phi  \|_{\mathcal{L}_2 (V) }  < \infty, \quad  \text{ for some } \beta \in (1 - 2H, 1].
\end{equation}
%
\end{assumption}
It is worthwhile to point out that the above series \eqref{eq:Wiener-representation} may not converge in $V$, 
but in some space $\tilde{V}$ into which $V$ can be embedded \cite{da2014stochastic}. 

\begin{assumption}[Initial value] \label{ass:X0}
Let $X_0 \colon  \Omega \rightarrow V $ be a $\mathcal{F}_0/\mathcal{B}(V)$-measurable mapping with  
$X_0 \in L^2 (\Omega, V_{ 2 H + \beta - 1 } )$.
\end{assumption}
We remark that the requirement of smooth initial data is not essential and one can reduce it at the expense of having the constant $C$ later depending on $T^{-1}$, by exploring the smoothing effect of the semigroup $E(t), t \in [0, T ]$ and standard nonsmooth data estimates. In this paper we prefer the smooth initial data to simplify the presentation, also to obtain constants uniform with respect to $T$ .
The above setting suffices to establish the following regularity results.
\begin{thm} \label{thm:SPDE-regularity-result}
Under Assumptions \ref{ass:A}-\ref{ass:X0}, SPDE \eqref{eq:spde} possesses a unique mild solution 
determined by \eqref{eq:intro-mild}. Furthermore, we have
\begin{equation} \label{eq:nonlinear-space-regularity}
\|  X(t) \|_{L^2 ( \Omega; V_{2H + \beta - 1}) } \leq C ( 1 + \| \xi \|_{L^2 ( \Omega; V_{2H + \beta - 1}) }),
\quad \text{ for } t \in [0, T]
\end{equation}
and for $\delta \in [0, 2H + \beta - 1 ]$, $t \geq s $,
\begin{equation}\label{eq:SPDE-optimal-time-regularity}
\|  X(t) - X(s) \|_{L^2 ( \Omega; V_{\delta}) } \leq C ( 1 + \| \xi \|_{L^2 ( \Omega; V_{2H + \beta - 1}) }) (t -s)^{ \frac{ 2H + \beta - 1 - \delta }{2} }.
\end{equation}
\end{thm}
The proof of Theorem \ref{thm:SPDE-regularity-result} is postponed to the end of this section.
Before that, we present an important lemma, which plays an essential role in 
deriving the sharp regularity results.
\begin{lem} \label{lem:optimal-regularity-lemma}
Under Assumption \ref{ass:A},  there exists a constant $C$ only depending on $H$ such that,
for $\delta \in [0, H]$, $0< s < t $ and $x \in V$, 
\begin{equation} \label{eq:optimal-lemma}
\int_s^t \int_s^t \langle A^{\delta} E( t - u ) x ,  A^{\delta} E ( t - v ) x \rangle \phi( u - v ) \, \dd u \dd v  
\leq
C ( t - s )^{ 2 (H - \delta) } \| x \|^2.
\end{equation}
\end{lem}

{\it Proof of Lemma \ref{lem:optimal-regularity-lemma}.} 
Due to the expansion of $x \in V$ in terms of the eigenbasis $\{ e_i \}_{i \in \N } $  of the operator $A$, 
one can arrive at
\begin{equation}\label{eq:optimal-proof1}
\begin{split}
& \int_s^t \int_s^t \langle A^{\delta} E( t - u ) x ,  A^{\delta} E ( t - v ) x \rangle \phi( u - v ) \, \dd u \dd v 
\\ 
& \quad =
\int_s^t \int_s^t \bigg \langle  \sum_{i \in \N}  \lambda_i ^{\delta} e^{ - \lambda_i ( t - u )} \langle x, e_i \rangle e_i ,  
\sum_{j \in \N}  \lambda_j ^{\delta} e^{ - \lambda_j ( t - v )} \langle x, e_j \rangle  e_j \bigg \rangle \phi( u - v ) \, \dd u \dd v 
\\
& \quad =
\int_s^t \int_s^t   \sum_{i \in \N}  \lambda_i ^{ 2 \delta } e^{ - \lambda_i ( 2 t - u - v ) } \langle x, e_i \rangle^2  \phi( u - v ) \, \dd u \dd v 
\\
& \quad =
\alpha_H \sum_{i \in \N}  \lambda_i ^{ 2 \delta }  \langle x, e_i \rangle^2 \int_s^t \int_s^t  e^{ - \lambda_i ( 2 t - u - v ) }  | u - v |^{ 2 H - 2} \, \dd u \dd v 
\\
& \quad =
\alpha_H \sum_{i \in \N}  \lambda_i ^{ 2 \delta }  \langle x, e_i \rangle^2 \int_s^t \int_s^v  e^{ - \lambda_i ( 2 t - u - v ) }  ( v - u )^{ 2 H - 2} \, \dd u \dd v
\\
& \qquad +
\alpha_H \sum_{i \in \N}  \lambda_i ^{ 2 \delta }  \langle x, e_i \rangle^2 \int_s^t \int_v^t  e^{ - \lambda_i ( 2 t - u - v ) }  ( u - v )^{ 2 H - 2} \, \dd u \dd v
\\
& \quad : =
I_1 + I_2.
\end{split}
\end{equation}
By a change of variable $ u = v w + s ( 1 - w ) $ we start the estimate of $I_1$:
\begin{equation}
I_ 1 = \alpha_H \sum_{i \in \N} \langle x, e_i \rangle^2 \lambda_i ^{ 2 \delta }  \int_s^t \int_0^1 
e^{ - \lambda_i [ 2 t - v w - s ( 1 - w ) - v ] } ( 1 - w )^{ 2 H - 2} ( v - s )^{ 2 H - 1} 
\, \dd w \dd v.
\end{equation}
Letting $ v = t r + s ( 1 - r ) $ further shows
\begin{equation}
\begin{split}
I_1 & =  \alpha_H \sum_{i \in \N} \langle x, e_i \rangle^2 \lambda_i ^{ 2 \delta }  \int_0^1 \int_0^1 e^{ - \lambda_i ( 2 - r w - r ) ( t -s ) } ( 1 - w )^{ 2 H - 2}  r ^{ 2 H - 1} ( t - s )^{ 2 H } \, \dd w \dd r
\\
& =
\alpha_H \sum_{ \lambda_i (t -s) \leq 1} \langle x, e_i \rangle^2 ( t - s )^{ 2 H }  \lambda_i ^{ 2 \delta }  \int_0^1 \int_0^1 e^{ - \lambda_i ( 2 - r w - r ) ( t -s ) } ( 1 - w )^{ 2 H - 2}  r ^{ 2 H - 1}  \, \dd w \dd r 
\\
& \quad + 
\alpha_H \sum_{ \lambda_i (t -s) > 1} \langle x, e_i \rangle^2 ( t - s )^{ 2 H }  \lambda_i ^{ 2 \delta }  \int_0^1 \int_0^1 e^{ - \lambda_i ( 2 - r w - r ) ( t -s ) } ( 1 - w )^{ 2 H - 2}  r ^{ 2 H - 1}  \, \dd w \dd r .
\\
& : = I_{11} + I_{12}.
\end{split}
\end{equation}
%
The estimate of the first term $I_{11}$ is quite easy:
\begin{equation} \label{eq:opt-lem-I11}
\begin{split}
I_{11} & \leq  \alpha_H \sum_{ \lambda_i (t -s) \leq 1} \langle x, e_i \rangle^2 ( t - s )^{ 2 ( H - \delta ) }  \int_0^1 \int_0^1 e^{ - \lambda_i ( 2 - r w - r ) ( t -s ) } ( 1 - w )^{ 2 H - 2}  r ^{ 2 H - 1}  \, \dd w \dd r 
\\
& \leq 
\alpha_H ( t - s )^{ 2 ( H - \delta ) }  \sum_{ \lambda_i (t -s) \leq 1} \langle x, e_i \rangle^2   \int_0^1 \int_0^1  ( 1 - w )^{ 2 H - 2}  r ^{ 2 H - 1}  \, \dd w \dd r 
\\
&=
\tfrac12 ( t - s )^{ 2 ( H - \delta ) }  \sum_{ \lambda_i (t -s) \leq 1} \langle x, e_i \rangle^2.
\end{split}
\end{equation}
Subsequently we handle the estimate of $I_{12}$. By exploiting elementary arguments such as integration by parts we have
\begin{equation}
\small
\begin{split}
I_{12} & \leq  \alpha_H \sum_{ \lambda_i (t -s) > 1} \langle x, e_i \rangle^2 ( t - s )^{ 2 H }  \lambda_i ^{ 2 \delta }   \int_0^1 e^{ - \lambda_i ( 1 - r ) ( t -s ) } r ^{ 2 H - 1} \dd r \int_0^1 e^{ - \lambda_i ( 1 -  w ) ( t -s ) } ( 1 - w )^{ 2 H - 2}   \, \dd w 
\\
& =
\alpha_H \sum_{ \lambda_i (t -s) > 1} \langle x, e_i \rangle^2 ( t - s )^{ 2 ( H - \delta ) }  [ \lambda_i ( t -s ) ]^{ 2 \delta - 1 } 
\Big[
1 - ( 2H -1 ) \int_0^1 e^{ - \lambda_i ( 1 - r ) ( t -s )  } r ^{ 2 H - 2 } \dd r
\Big]
\\ & \qquad 
\times
\int_0^1 e^{ - \lambda_i ( t -s ) w } w^{ 2 H - 2}   \, \dd w
\\ 
& \leq
\alpha_H \sum_{ \lambda_i (t -s) > 1} \langle x, e_i \rangle^2 ( t - s )^{ 2 ( H - \delta ) }  [ \lambda_i ( t -s ) ]^{ 2 \delta - 1 } 
\int_0^1 e^{ - \lambda_i ( t -s ) w } w^{ 2 H - 2}   \, \dd w
\\ 
& =
\alpha_H \sum_{ \lambda_i (t -s) > 1} \langle x, e_i \rangle^2 ( t - s )^{ 2 ( H - \delta ) }  [ \lambda_i ( t -s ) ]^{ 2 ( \delta - H) } 
\int_0^{ \lambda_i ( t -s ) } e^{ - \theta }  \theta^{ 2 H - 2}   \, \dd \theta
\\ 
& \leq
\alpha_H ( t - s )^{ 2 ( H - \delta ) }  \sum_{ \lambda_i (t -s) > 1} \langle x, e_i \rangle^2  
\int_0^{ \infty } e^{ - \theta }  \theta^{ 2 H - 2}   \, \dd \theta
\\
& \leq
c_H ( t - s )^{ 2 ( H - \delta ) }  \sum_{ \lambda_i (t -s) > 1} \langle x, e_i \rangle^2.
\end{split}
\end{equation}
Combining the above two estimates together yields
\begin{equation} \label{eq:I1-final-estimate}
I_1 \leq I_{11} + I_{12} \leq c_H ( t - s )^{ 2 ( H - \delta ) }  \sum_{ i \in \N} \langle x, e_i \rangle ^2
=
c_H ( t - s )^{ 2 ( H - \delta ) }  \| x \|^2.
\end{equation}
With regard to $I_2$, one can observe that, after an alternative representation of integral domain,
\begin{equation}\label{eq:Optimal-regularity-I2-estimate}
I_2 = \alpha_H \sum_{i \in \N}  \lambda_i ^{ 2 \delta }  \langle x, e_i \rangle^2 \int_s^t \int_s^u  e^{ - \lambda_i ( 2 t - u - v ) }  ( u - v )^{ 2 H - 2} \, \dd v \dd u = I_1,
\end{equation}
which in conjunction with \eqref{eq:I1-final-estimate} implies the assertion as required. $\square$

\begin{prop} \label{prop:optimal-space-regularity-stoch-conv}
Under Assumptions \ref{ass:A} and \ref{ass:Phi} with $\beta \in ( 1 - 2H, 1]$, the stochastic convolution
\begin{equation}
\mathcal{O}_t := 
\int_0^t E(t-s) \Phi \, \dd W^H(s), \quad t \in [0, T]
\end{equation}
is well-defined in $L^2 ( \Omega; V)$ for any $t \in [0, T]$ and possesses the following spatial regularity:
\begin{equation}\label{eq:optimal-space-regularity-stoch-conv}
\|  \mathcal{O}_t \|_{L^2 ( \Omega; V_{2H + \beta - 1}) } 
\leq 
c_H \| A^{ \frac{\beta -1}{2} } \Phi  \|_{\mathcal{L}_2 (V) }, 
\quad
\forall \: t \in [0, T].
\end{equation}
Furthermore, for $ 0 \leq s < t \leq T$ the following temporal regularity holds:
\begin{equation}\label{eq:optimal-time-regularity-stoch-conv}
\|  \mathcal{O}_t - \mathcal{O}_s \|_{L^2 ( \Omega; V_{\delta}) } \leq c_H \| A^{ \frac{\beta -1}{2} } \Phi  \| _{\mathcal{L}_2 (V) }
(t -s)^{ \frac{ 2H + \beta - 1 - \delta }{2} },
\quad
\forall \, \delta \in [0, 2H + \beta - 1 ].
\end{equation}
\end{prop}
{\it Proof of Proposition \ref{prop:optimal-space-regularity-stoch-conv}. }
The well-posedness of the stochastic convolution $\mathcal{O}_t$, $t \in [0, T]$ is clear 
by taking the knowledges in section \ref{sec:prelim} into account. This has also been confirmed in \cite{duncan2002fractional}. So we just check \eqref{eq:optimal-space-regularity-stoch-conv}
and \eqref{eq:optimal-time-regularity-stoch-conv}.
Thanks to the It\^o isometry \eqref{eq:intro-ito-isometry} and Lemma \ref{lem:optimal-regularity-lemma},
\begin{equation} \label{eq:proof-optimal-space-stoch-conv}
\small
\begin{split}
\|  \mathcal{O}_t \|_{L^2 ( \Omega; V_{2H + \beta - 1}) }^2 
& =
\int_0^t \int_0^t \left\langle A^{ \frac{2H + \beta - 1}{2} } E( t - u )  \Phi  ,  
A^{ \frac{2H + \beta - 1}{2} }  E ( t - v )  \Phi \right\rangle_{ \mathcal{L}_2 } \phi( u - v ) \, \dd u \dd v 
\\
& = 
\sum_{i \in \N} 
\int_0^t \int_0^t \left\langle A^H E( t - u ) A^{ \frac{\beta - 1}{2} }  \Phi  e_i ,  
A^H  E ( t - v ) A^{ \frac{\beta - 1}{2} }  \Phi e_i \right\rangle \phi( u - v ) \, \dd u \dd v 
\\
& \leq
c_H \sum_{i \in \N}  \| A^{ \frac{\beta - 1}{2} }  \Phi e_i \|^2 
\\
& = c_H \| A^{ \frac{\beta -1}{2} }  \Phi  \|_{\mathcal{L}_2}^2,
\end{split}
\end{equation}
where the fact was also used that densely defined linear operators commute 
with the stochastic integral \cite[Proposition 2.4]{duncan2009semilinear}.
%
%
To get \eqref{eq:optimal-time-regularity-stoch-conv},  we first realize that
\begin{equation}
\mathcal{O}_t - \mathcal{O}_s = \big( E( t -s ) - I \big) \mathcal{O}_s + \int_s^t E( t - r )  \Phi \, \dd W^H(r).
\end{equation}
Using \eqref{eq:E.Inequality} and \eqref{eq:proof-optimal-space-stoch-conv} shows
\begin{equation}
\begin{split}
\|  \big( E( t -s ) - I \big) \mathcal{O}_s \|_{L^2 ( \Omega; V_{\delta}) }
& \leq 
\|  \big( E( t -s ) - I \big) A^{ \frac{\delta - ( 2H + \beta - 1)}{2} } \|_{\mathcal{L}(V) } 
\cdot
\|  \mathcal{O}_s \|_{L^2 ( \Omega; V_{2H + \beta - 1}) }
\\ 
& \leq
C ( t - s )^{ \frac{ 2H + \beta - 1 - \delta }{2} } \| A^{ \frac{\beta -1}{2} }  \Phi \|_{\mathcal{L}_2}.
\end{split}
\end{equation}
Again, It\^o's isometry \eqref{eq:intro-ito-isometry} and Lemma \ref{lem:optimal-regularity-lemma} imply that
\begin{equation}\label{eq:Ot-Holder-estimate2}
\begin{split}
& \Big\|  \int_s^t E( t - r )  \Phi  \, \dd W^H(r) \Big\|_{L^2 ( \Omega; V_{\delta}) }^2
 =
\Big\|  \int_s^t A^{\frac{\delta}{2} } E( t - r )  \Phi \, \dd W^H(r) \Big\|_{L^2 ( \Omega; V) }^2
\\
& \quad =
\int_s^t \int_s^t \left\langle A^{ \frac{\delta}{2} } E( t - u )  \Phi,  
A^{ \frac{\delta}{2} }  E ( t - v )  \Phi \right\rangle_{ \mathcal{L}_2 } \! \phi( u - v ) \, \dd u \dd v 
\\
& \quad =
\sum_{ i \in \N }
\int_s^t \int_s^t \left\langle A^{ \frac{\delta + 1 - \beta}{2} } E( t - u ) A^{ \frac{\beta -1}{2} }  \Phi  e_i ,  
A^{ \frac{\delta + 1 - \beta}{2} }  E ( t - v ) A^{ \frac{\beta -1}{2} }  \Phi e_i \right\rangle \! \phi( u - v ) \, \dd u \dd v 
\\ 
& \quad \leq
C ( t - s )^{ 2H + \beta - 1 - \delta } \| A^{ \frac{\beta -1}{2} }  \Phi  \|_{\mathcal{L}_2}^2.
\end{split}
\end{equation}
The assertion \eqref{eq:optimal-time-regularity-stoch-conv} straightforwardly follows from 
the triangle inequality. $\square$
\begin{rem}
Before moving on, we present another simple but crude way to analyze the spatial regularity of $ \mathcal{O}_t $.
Using It\^o's isometry, some basic inequalities \eqref{eq:L-L2-norm-inequality} 
and Lemma \ref{lem:double-intefral-inequality} promises
\begin{equation*}
\begin{split}
\|  \mathcal{O}_t \|_{L^2 ( \Omega; V_{\rho}) }^2
& =
\int_0^t \int_0^t \left\langle A^{ \frac{\rho}{2} } E( t - u )  \Phi  ,  
A^{ \frac{\rho}{2} }  E ( t - v )  \Phi \right\rangle_{ \mathcal{L}_2 } \phi( u - v ) \, \dd u \dd v 
\\
& \leq
\int_0^t \int_0^t \big\| A^{ \frac{\rho + 1 - \beta}{2} } E( t - u )  A^{ \frac{\beta - 1}{2} }  \Phi \big\|_{ \mathcal{L}_2 }  
 \big\| A^{ \frac{\rho + 1 - \beta}{2} }  E ( t - v )  A^{ \frac{\beta - 1}{2} }  \Phi \big\|_{ \mathcal{L}_2 } \phi( u - v ) \, \dd u \dd v 
\\
& \leq
\|  A^{ \frac{\beta - 1}{2} }  \Phi \|^2_{ \mathcal{L}_2 }
\int_0^t \int_0^t \big\| A^{ \frac{\rho + 1 - \beta}{2} } E( t - u )  \big\|_{ \mathcal{L}(V) }  
\big\| A^{ \frac{\rho + 1 - \beta}{2} }  E ( t - v ) \big\|_{ \mathcal{L}(V) }  
\phi( u - v ) \, \dd u \dd v  
\\
& \leq
C \|  A^{ \frac{\beta - 1}{2} }  \Phi \|^2_{ \mathcal{L}_2 }
\int_0^t \int_0^t 
(t - u)^{- \frac{\rho + 1 - \beta}{2} } 
(t - v)^{- \frac{\rho + 1 - \beta}{2} }  
\phi( u - v ) \, \dd u \dd v  
<
\infty
\end{split}
\end{equation*}
for any $\rho < 2 H + \beta - 1$.
Apparently, this way leads us to reduced spatial mean-square regularity for the stochastic convolution $ \mathcal{O}_t $, with the border case $\rho = 2 H + \beta - 1$ not included.
\end{rem}

We are now in a position to verify Theorem \ref{thm:SPDE-regularity-result}.

{\it Proof of Theorem \ref{thm:SPDE-regularity-result}.}  
Bearing in mind that $ \mathcal{O}_t, t \in [0, T] $ is well-defined in 
$L^2 ( \Omega; V)$ and that the nonlinear mapping $F$ obeys the globally Lipschitz
condition, one can follow standard arguments to acquire the existence and 
the uniqueness of a mild solution to \eqref{eq:spde} in $L^2 ( \Omega; V)$, given by \eqref{eq:intro-mild}. 
So we just focus on the regularity results of the mild solution.
By \eqref{eq:E.Inequality} and \eqref{eq:optimal-space-regularity-stoch-conv},
\begin{equation}
\begin{split}
\|  X(t) \|_{L^2 ( \Omega; V_{2H + \beta - 1}) } 
& \leq
\| E(t) \xi \|_{L^2 ( \Omega; V_{2H + \beta - 1}) }
    +
   \Big \| \int_0^t E(t-s) F ( X( s ) ) \, \dd s \Big \|_{L^2 ( \Omega; V_{2H + \beta - 1}) }
\\ 
   & \quad +
    \|  \mathcal{O}_t \|_{L^2 ( \Omega; V_{2H + \beta - 1}) }
\\ & \leq
\| \xi \|_{L^2 ( \Omega; V_{2H + \beta - 1}) }
+
C \int_0^t (t-s)^{ -\frac{2H + \beta - 1}{2} } \big \| F ( X( s ) ) \big \|_{L^2 ( \Omega; V) } \, \dd s 
\\
  & \quad
+
c_H \| A^{ \frac{\beta -1}{2} }  \Phi  \|_{\mathcal{L}_2}
\\ & \leq
C ( 1 + \| \xi \|_{L^2 ( \Omega; V_{2H + \beta - 1}) }).
\end{split}
\end{equation}
This confirms \eqref{eq:nonlinear-space-regularity}.  Concerning the proof of \eqref{eq:SPDE-optimal-time-regularity},
one can easily see that
\begin{equation}
\begin{split}
X(t) - X(s)  
& =
\big( E(t - s) - I \big) X(s) 
+
\int_s^t \! E(t-r) F ( X( r ) ) \, \dd r
+
\int_s^t \! E( t - r )  \Phi \, \dd W^H(r),
\end{split}
\end{equation}
and thus we used \eqref{eq:Ot-Holder-estimate2} to obtain
\begin{equation}
\begin{split}
\|  X(t) - X(s) \|_{L^2 ( \Omega; V_{\delta}) } 
& \leq 
\big \|  \big(  E(t - s) - I \big) X(s) \big \|_{L^2 ( \Omega; V_{\delta}) }
\\
& \quad +
\Big \| \int_s^t E(t-r) F ( X( r ) ) \, \dd r \Big \|_{L^2 ( \Omega; V_{\delta}) }
\\
& \quad +
\Big \| \int_s^t E( t - r )  \Phi \, \dd W^H(t) \Big \|_{L^2 ( \Omega; V_{\delta}) }
\\
& \leq
\big \|  \big(  E(t - s) - I \big) A^{ \frac{ \delta - 2H - \beta + 1 }{2} } \big \|_{\mathcal{L} (V) } 
\big \| X(s) \big \|_{L^2 ( \Omega; V_{2H + \beta - 1}) }
\\
& \quad +
C ( t - s)
+
C ( t - s )^{ \frac{ 2H + \beta - 1 - \delta }{2} } \| A^{ \frac{\beta -1}{2} }  \Phi  \|_{\mathcal{L}_2 }
\\
& \leq
C ( 1 + \| \xi \|_{L^2 ( \Omega; V_{2H + \beta - 1}) })  (t -s)^{ \frac{ 2H + \beta - 1 - \delta }{2} },
\end{split}
\end{equation}
which completes the proof of Theorem \ref{thm:SPDE-regularity-result}. $\square$

\section{Optimal convergence rates of a full-discretization}
\label{sec:strong.convergence}
This section is devoted to the error analysis for strong approximations of the underlying problem.
Optimal convergence rates are obtained in the mean-square sense, which coincide with previously 
derived regularity of the mild solution.

\subsection{Spatial semi-discretization}
In this part we spatially discretize \eqref{eq:spde}  with a spectral Galerkin method.
To this end, for $N\in \mathbb{N}$ we define a finite dimensional subspace of $V$ by
$
 V^N := \mbox{span} \{e_1, e_2, \cdots, e_N \},
$
and the projection operator $P_N \colon V_{\alpha}\rightarrow V^N$ by
\begin{align}
  P_N \xi = \sum_{i=1}^N \langle \xi, e_i \rangle e_i, \quad \forall\, \xi \in V_{\alpha}, \, \alpha \in \R.
\end{align}
Here $V^N$ is chosen as the linear space spanned by the $N$ first eigenvectors of 
the dominant linear operator $A$. It is easy to see that
\begin{equation}\label{eq:P-N-estimate}
\| ( P_N - I ) \varphi \|  \leq \lambda_{N+1}^{- \frac {\alpha}{2} } \|\varphi\|_{\alpha}, 
\quad \forall \: \varphi \in V_{\alpha}, \: \alpha \geq0.
\end{equation}
Additionally, define $A_N \colon V \rightarrow V^N$ as
$A_N = A P_N $, which generates an analytic semigroup $E_N(t) = e^{-t A_N}$, $t \in [0, \infty)$ in $V^N$. 
Then the spectral Galerkin method for \eqref{eq:spde} is described by
\begin{equation}\label{eq:spectral-spde}
\begin{split}
\left\{
    \begin{array}{lll} \dd X^N(t) + A_N X^N ( t ) \dd t = P_N F ( X^N (t) ) \dd t  +  P_N \Phi \, \dd W^H (t), 
    \quad t \in (0, T], \\
     X^N(0) = P_N \xi.
    \end{array}\right.
\end{split}
\end{equation}
The corresponding mild solution is given by
\begin{equation}\label{eq:Spectral-Galerkin-mild}
    X^N(t)
    =
    E_N(t) P_N \xi
    +
    \int_0^t E_N (t - s ) P_N F ( X^N( s ) ) \, \dd s
    +
    \int_0^t E_N(t-s) P_N\Phi \, \dd W^H(s), \:  \mathbb{P} \mbox{-a.s.}.
\end{equation}
Noting that
\begin{equation}\label{eq:EN-P-N-commun}
E_N(t) P_N = E (t) P_N,
\end{equation}
and taking \eqref{eq:E.Inequality} and \eqref{eq:P-N-estimate} into account promise 
\begin{equation}\label{eq:spatial-operator-error}
\| ( E(t) -  E_N(t) P_N )  \varphi \| \leq C_{\alpha} t^{ - \frac{\alpha}{2} } 
\lambda_{N+1}^{- \frac {\alpha}{2} } \| \varphi \|,
\quad
t > 0, \, \alpha \geq 0, \, \varphi \in V.
\end{equation}
The error analysis of the spatial discretization \eqref{eq:spectral-spde}
also relies on the following error estimate.
\begin{lem}
\label{lem:Spectral-Galerkin-lemma}
For $\beta \in ( 1 - 2H, 1]$, it holds that
\begin{equation} \label{eq:optimal-Spectral-Galerkin-lemma}
\int_0^t \int_0^t  \Big \langle  \mathcal{K}_N ( t - u ) x ,   \mathcal{K}_N ( t - v ) x  \Big \rangle \phi( u - v ) \, \dd u \dd v  
\leq
C  \lambda_{N + 1}^{- ( 2 H + \beta - 1) } \| x \|_{\beta - 1}^2,
\end{equation}
where the error operator $\mathcal{K}_N ( t )$ is defined by $\mathcal{K}_N ( t ) x : = \big( E ( t ) - E_N ( t ) P_N \big) x$, for $x \in V_{\beta - 1} $.
\end{lem}
{\it Proof of Lemma \ref{lem:Spectral-Galerkin-lemma}.}
By virtue of Lemma \ref{lem:optimal-regularity-lemma} with $\delta = H$, one can achieve
\begin{equation}
\begin{split}
& \int_0^t \int_0^t \langle  \mathcal{K}_N ( t - u ) x ,   \mathcal{K}_N ( t - v ) x \rangle \phi( u - v ) \, \dd u \dd v 
\\ 
& \quad =
\int_0^t \int_0^t \bigg \langle  \sum_{i \in \N}  \langle \mathcal{K}_N ( t - u ) x, e_i \rangle  e_i ,  
\sum_{j \in \N}   \langle \mathcal{K}_N ( t - v ) x, e_j \rangle  e_j \bigg \rangle \phi( u - v ) \, \dd u \dd v 
\\
& \quad =
\int_0^t \int_0^t   \sum_{i \in \N}  \langle \mathcal{K}_N ( t - u ) x, e_i \rangle 
\langle \mathcal{K}_N ( t - v ) x, e_i \rangle \phi( u - v ) \, \dd u \dd v
\\
& \quad =
\int_0^t \int_0^t   \sum_{i \geq N+1} e^{- \lambda_i ( 2 t - u - v) } \langle x, e_i \rangle^2  \phi( u - v ) \, \dd u \dd v 
\\
& \quad \leq
\lambda_{N + 1}^{- ( 2 H + \beta - 1) }
\int_0^t \int_0^t   \sum_{i \geq N+1} \lambda_{i}^{ 2 H + \beta - 1 } e^{- \lambda_i ( 2 t - u - v) } 
\langle x, e_i \rangle^2  \phi( u - v ) \, \dd u \dd v
\\
& \quad \leq
\lambda_{N + 1}^{- ( 2 H + \beta - 1) }
\int_0^t \int_0^t   \sum_{i \in \N } \lambda_{i}^{ 2 H } e^{- \lambda_i ( 2 t - u - v) } 
\langle A^{\frac{\beta - 1}{2}} x, e_i \rangle^2  \phi( u - v ) \, \dd u \dd v
\\
& \quad =
\lambda_{N + 1}^{- ( 2 H + \beta - 1) }
\int_0^t \int_0^t 
\Big \langle 
A^{H} E( t - u ) A^{\frac{\beta - 1}{2}} x ,  A^{H} E ( t - v ) A^{\frac{\beta - 1}{2}} x 
\Big \rangle \phi( u - v ) 
\, \dd u \dd v  
\\
& \quad \leq
C \lambda_{N + 1}^{- ( 2 H + \beta - 1) } \| x \|_{\beta - 1}.
\end{split}
\end{equation}
This validates \eqref{eq:optimal-Spectral-Galerkin-lemma}. $\square$

Equipped with the above lemma, we are now prepared to show the following convergence result for 
the spectral Galerkin discretization \eqref{eq:spectral-spde}.
\begin{thm}[Spatial error estimate]
\label{thm:space-main-conv}
Let Assumptions \ref{ass:A}-\ref{ass:X0} hold with $\beta \in ( 1 - 2H, 1]$. 
Let  $X(t)$ and $X^N(t)$ be defined through \eqref{eq:intro-mild} 
and \eqref{eq:Spectral-Galerkin-mild}, respectively. Then
\begin{equation}\label{eq:spatial-error}
\| X(t) - X^N(t) \|_{L^2 ( \Omega; V) } \leq
C ( 1 + \| \xi \|_{L^2 ( \Omega; V_{2H + \beta - 1}) })
\lambda_{N + 1}^{- \frac{ 2 H + \beta - 1}{2} }.
\end{equation}
\end{thm}
{\it Proof of Theorem \ref{thm:space-main-conv}. } Evidently,
\begin{equation}\label{eq:space-err-proof-split}
\begin{split}
   \| X(t) - X^N(t) \|_{L^2 ( \Omega; V) }
    & \leq
    \|  ( E(t) - E_N ( t ) P_N) \xi \|_{L^2 ( \Omega; V) }
    \\ & 
    \quad +
   \Big \|  \int_0^t  \Big(  E (t - s ) F ( X( s ) ) - E_N (t - s ) P_N F ( X^N( s ) ) \Big)  \dd s \Big \|_{L^2 ( \Omega; V) }
    \\ & 
    \quad +
     \Big \| \int_0^t \big( E(t-s)  - E_N( t - s ) P_N \big) \Phi \, \dd W^H(s) \Big \|_{L^2 ( \Omega; V) } 
    \\
    & := J_1 + J_2 + J_3.
\end{split}
\end{equation}
The estimate of $J_1$ is obvious after taking \eqref{eq:P-N-estimate} into account:
\begin{equation}\label{eq:J1}
J_1 \leq  \lambda_{N + 1}^{- \frac{ 2 H + \beta - 1}{2} } \| \xi \|_{L^2 ( \Omega; V_{2 H + \beta - 1} ) }.
\end{equation}
The second term $J_2$ can be treated as
\begin{equation}\label{eq:J2}
\begin{split}
J_2 & \leq  \Big \|  \int_0^t  \Big(  E (t - s )  - E_N (t - s ) P_N \Big) F ( X( s ) ) \dd s \Big \|_{L^2 ( \Omega; V) }
\\ & \quad +
\Big \|  \int_0^t   E_N (t - s ) P_N \Big( F ( X( s ) ) -  F ( X^N( s ) ) \Big)  \dd s \Big \|_{L^2 ( \Omega; V) }
\\ & \leq
C \lambda_{N + 1}^{- \frac{ 2 H + \beta - 1}{2} }   
\int_0^t (t -s)^{- \frac{2H + \beta -1}{2 } } \| F (X (s) ) \|_{L^2 ( \Omega; V) } \,\dd s
+
L \int_0^t  \big \|  X( s ) -  X^N( s ) \big \|_{L^2 ( \Omega; V) } \dd s
\\ & \leq
C \big( 1 + \| \xi \|_{L^2 ( \Omega; V) } \big) \lambda_{N + 1}^{- \frac{ 2 H + \beta - 1}{2} } 
+ 
L \int_0^t  \big \|  X( s ) -  X^N( s ) \big \|_{L^2 ( \Omega; V) } \dd s,
\end{split}
\end{equation}
where Assumption \ref{ass:F}, Theorem \ref{thm:SPDE-regularity-result}  and 
\eqref{eq:spatial-operator-error} were employed.
It remains to deal with $J_3$. To do so we use  It\^o's isometry \eqref{eq:intro-ito-isometry},  
and Lemma \ref{lem:Spectral-Galerkin-lemma} to arrive at
\begin{equation}\label{eq:J3}
\begin{split}
| J_3 |^2  & =  \int_0^t \int_0^t \Big \langle   \mathcal{K}_N ( t - u ) \Phi,  
  \mathcal{K}_N ( t - v ) \Phi  \Big \rangle_{ \mathcal{L}_2 } \! \phi( u - v ) \, \dd u \dd v 
\\  & =
\sum_{ i \in \N} \int_0^t \int_0^t \big \langle   \mathcal{K}_N ( t - u ) \Phi e_i,  
  \mathcal{K}_N ( t - v ) \Phi e_i \big \rangle  \phi( u - v ) \, \dd u \dd v 
\\ & \leq
C  \lambda_{N + 1}^{- ( 2 H + \beta - 1) } \| A^{\frac{\beta - 1}{2} } \Phi \|_{\mathcal{L}_2 }^2.
\end{split}
\end{equation} 
Inserting \eqref{eq:J1}, \eqref{eq:J2} and \eqref{eq:J3} into \eqref{eq:space-err-proof-split} 
and invoking the Gronwall inequality give \eqref{eq:spatial-error}. 
$\square$
\subsection{A full discretization}
This subsection concerns a time discretization of the spatially discretized problem \eqref{eq:spectral-spde}.
For $M \in \N$ we construct a uniform mesh on $[0, T]$ with $\tau = \tfrac{T}{M}$ being the stepsize.
Applying the linear implicit Euler time discretization to \eqref{eq:spectral-spde} 
results in a spatio-temporal full discretization:
\begin{equation}\label{eq:full-discrete-scheme-Spectral}
\bar{X}^N_{m+1} = R (\tau A_N) \bar{X}^N_{m} + \tau R (\tau A_N) P_N F ( \bar{X}^N_{m} ) 
+ R (\tau A_N) P_N  \Phi \Delta W^H_m, \ m = 0, 1,..., M-1,
\end{equation}
where $\Delta W^H_m : = W^H ( t_{m+1} ) -  W^H ( t_{m} ) $ are the Brownian increments
and define $R (z) := ( 1 + z )^{-1}, z \geq 0$.
As shown in the proof of \cite[Theorem 7.1]{thomee2006galerkin},  there exists a constant $C \geq 0$ such that
\begin{equation}
| R (z) - e^{ - z } | \leq C z^{2}, \quad z \in [0, 1],
\end{equation}
and there exists a constant $c \in (0, 1)$ such that
\begin{equation}\label{eq:Rz-exponential}
| R (z) | \leq e^{ - c z }, \quad z \in [0, 1].
\end{equation}
These two inequalties suffice to ensure that, for $j =1,2,3,...$,
\begin{equation} \label{eq:rational-appro-error}
| R ( z )^j - e^{- z j } | 
= 
\Big| \big( R ( z ) - e^{- z } \big) 
\sum_{i = 0}^{ j -1 }  R ( z )^{j - 1 - i} e^{ - z i }  \Big |
\leq
C j z ^2 e^{ - c ( j - 1 ) z },
\quad z \in [0, 1],
\end{equation}
where $c \in (0, 1)$ comes from \eqref{eq:Rz-exponential}.
Employing the above facts one can show that
\begin{equation} \label{eq:rational-exp-differ}
\big\| \big( R ( \tau A )^m - E (t_m) \big) x \big\| 
\leq 
C \tau^{ \frac{\mu}{2} } t_m^{- \frac{\mu-\nu}{2} } \| x \|_{\nu},
\quad
0 \leq \nu \leq \mu \leq 2,
x \in V_{\nu}.
\end{equation}
At the moment we state the error estimate for the time-stepping scheme \eqref{eq:full-discrete-scheme-Spectral} 
and its proof is put after Lemma \ref{lem:optimal-regularity-X-N}.
\begin{thm}[Convergence rates of temporal discretization]
\label{thm:time-convergence}
Let Assumptions \ref{ass:A}-\ref{ass:X0} hold with $\beta \in ( 1 - 2H, 1]$. Let $X^N (t_m) $ and $\bar{ X}^N_{m} $ 
be given by \eqref{eq:Spectral-Galerkin-mild} and \eqref{eq:full-discrete-scheme-Spectral}, respectively. 
Then
\begin{equation}\label{eq:thm-time-error}
\| X^N ( t_{m} ) -  \bar{ X}^N_{m} \|_{L^2 ( \Omega; V) }
     \leq
C \big( 1 + \| \xi \|_{L^2 ( \Omega; V_{2H + \beta - 1}) } \big) \tau^{ \frac{ 2 H + \beta - 1 }{2} }.
\end{equation}
\end{thm}
A combiniation of this with Theorem \ref{thm:space-main-conv} implies error bounds for the full discretization.
\begin{cor}[Error bounds for full discretization]
\label{cor:full-discrete-scheme-error-bound}
Under Assumptions \ref{ass:A}-\ref{ass:X0}, it holds that
\begin{equation}
\| X(t_m) - \bar{ X}^N_{m} \|_{L^2 ( \Omega; V) } \leq
C \big( 1 + \| \xi \|_{L^2 ( \Omega; V_{2H + \beta - 1}) } \big) 
\big(  \lambda_{N + 1}^{- \frac{ 2 H + \beta - 1}{2} } +
\tau^{ \frac{ 2 H + \beta - 1 }{2} } \big).
\end{equation}
\end{cor}
To launch the proof of Theorem \ref{thm:time-convergence}, we require the following ingredients. 
%
\begin{lem} \label{lem:phi-estimate}
For $i, j \in \N$, it holds that
\begin{equation}
\int_0^1 \int_0^1 \phi ( u + i - v - j ) \, \dd u \dd v  = 1 \quad \text{ for } \: i = j,
\end{equation}
and 
\begin{equation}\label{eq:lem-phi-estimate}
\int_0^1 \int_0^1 \phi ( u + i - v - j ) \, \dd u \dd v  \leq  \tfrac12 \max(i, j)^{2H - 1} \quad \text{ for } \: i \neq j.
\end{equation}
\end{lem}
{\it Proof of Lemma \ref{lem:phi-estimate}. } 
Recall that $\phi( y ) = \alpha_H | y |^{2 H - 2}, y \in \R$, with $\alpha_H = H ( 2 H - 1)$, $H \in (\tfrac12, 1)$. 
For the case $i = j$, elementary calculations yield
\begin{equation}
\small
\begin{split}
\int_0^1 \int_0^1 \phi ( u + i - v - j ) \, \dd u \dd v 
& = 
\alpha_H  \int_0^1 \int_0^v (v  - u )^{2H - 2} \, \dd u \dd v
+
\alpha_H  \int_0^1 \int_v^1 (u  - v )^{2H - 2} \, \dd u \dd v
=
1.
\end{split}
\end{equation}
For the case $i \neq j$, we assume $i > j$ without loss of generality.  
Then one can easily check that
\begin{equation}
\begin{split}
\int_0^1 \int_0^1 \phi ( u + i - v - j ) \, \dd u \dd v 
& = 
\alpha_H  \int_0^1 \int_0^1 (u + i - v - j )^{2H - 2} \, \dd u \dd v
\\ & =
\tfrac12 \big[ 
(i - j + 1)^{2H} - 2 ( i - j)^{2H} + ( i - j - 1)^{2H}
\big]
\\ & \leq
\tfrac12  (i - j + 1)^{2H-1},
\end{split}
\end{equation}
and \eqref{eq:lem-phi-estimate} is hence verified. $\square$
\begin{lem} \label{lem:lambda-phi-integral}
Let $\kappa_1, \kappa_2 \in \{ 0, 1\} $ and let $\lambda, t \in (0, \infty) $. 
Then there exists a uniform constant $C \in (0, \infty)$ independent of $\lambda, t $ such that
\begin{equation} \label{eq:lem-lambda-phi-integral}
\lambda^{ 2 H + \kappa_1 + \kappa_2 } \int_0^t \int_0^t
u^{ \kappa_1 }   v^{ \kappa_2 }  e^{ - \lambda ( u + v ) } \phi (u -v) \, \dd u \dd v
\leq C.
\end{equation}
\end{lem}
{\it Proof of Lemma \ref{lem:lambda-phi-integral}.} 
Note first that
\begin{equation}
\begin{split}
\int_0^t \int_0^t
u^{ \kappa_1 }   v^{ \kappa_2 }  e^{ - \lambda ( u + v ) } \phi (u -v) \, \dd u \dd v
& =  
\alpha_H \int_0^t \int_0^ v 
u^{ \kappa_1 }   v^{ \kappa_2 }  e^{ - \lambda ( u + v ) } (v - u)^{ 2 H - 2 } \, \dd u \dd v
\\ & \quad +
\alpha_H \int_0^t \int_v^ t
u^{ \kappa_1 }   v^{ \kappa_2 }  e^{ - \lambda ( u + v ) } (u -v)^{ 2 H - 2 } \, \dd u \dd v
\\ & := 
\mathbb{I}_1 + \mathbb{I}_2.
\end{split}
\end{equation}
A change of variable $ u = v r $ gives
\begin{equation} \label{eq:estimate-mathbb-I1}
\begin{split}
\mathbb{I}_1 & =  \alpha_H \int_0^t \int_0^ 1 
 v^{ 2H - 1 + \kappa_1 + \kappa_2 } r^{\kappa_1} e^{ - \lambda ( 1 + r ) v } (1 - r)^{ 2 H - 2 } \, \dd r \dd v
\\  & \leq
\alpha_H \int_0^t  
v^{ 2H - 1 + \kappa_1 + \kappa_2 }  e^{ - \lambda v }  \, \dd v 
\int_0^ 1  (1 - r)^{ 2 H - 2 } \, \dd r
\\ & =
H \lambda^{ - ( 2 H + \kappa_1 + \kappa_2 ) }
\int_0^{ \lambda t } 
\theta^{ 2H - 1 + \kappa_1 + \kappa_2 }  e^{ - \theta }  \, \dd \theta
\\ & \leq
C \lambda^{ - ( 2 H + \kappa_1 + \kappa_2 ) }.
\end{split}
\end{equation}
The same arguments as already used in \eqref{eq:Optimal-regularity-I2-estimate} and 
the estimate of  $\mathbb{I}_1$ above help us to obtain
\begin{equation} \label{eq:estimate-mathbb-I2}
\mathbb{I}_2
= 
\alpha_H \int_0^t \int_0^u
u^{ \kappa_1 }   v^{ \kappa_2 }  e^{ - \lambda ( u + v ) } (u -v)^{ 2 H - 2 } \, \dd v \dd u
\leq
C \lambda^{ - ( 2 H + \kappa_1 + \kappa_2 ) }.
\end{equation}
Gathering \eqref{eq:estimate-mathbb-I1} and \eqref{eq:estimate-mathbb-I2} finishes the proof of 
\eqref{eq:lem-lambda-phi-integral}.
$\square$

Repeating the proof of Theorem \ref{thm:SPDE-regularity-result}, one can show
\begin{lem}
\label{lem:optimal-regularity-X-N}
Let $X^N (t) $ be the solution to \eqref{eq:spectral-spde}, given by \eqref{eq:Spectral-Galerkin-mild}.
Under Assumptions \ref{ass:A}-\ref{ass:X0} with $\beta \in (1-2H, 1]$, it holds that, for $t, s \in [0, T]$ and $t > s$,
\begin{align} \label{eq:X-N-optimal-space}
\|  X^N(t) \|_{L^2 ( \Omega; V) } & \leq C \big( 1 + \| \xi \|_{L^2 ( \Omega; V) } \big),
\\
\label{eq:X-N-optimal-time}
\|  X^N(t) - X^N(s) \|_{L^2 ( \Omega; V) } & \leq  C \big( 1 + \| \xi \|_{L^2 ( \Omega; V_{2 H + \beta - 1}) } \big) 
(t -s)^{ \frac{ 2H + \beta - 1 }{2} }.
\end{align}
\end{lem}
We are now in a position to prove the main result, Theorem \ref{thm:time-convergence}.
\vspace{0.5cm}
\\
{\it Proof of Theorem \ref{thm:time-convergence}. }
Equivalently, the full discretization \eqref{eq:full-discrete-scheme-Spectral} can be expressed by
\begin{equation}
\bar{X}^N_{m} =  R (\tau A_N)^m P_N \xi + \tau \sum_{ i = 0}^{ m - 1}  R (\tau A_N)^{m-i} P_N F ( \bar{X}^N_{i} )
+ \sum_{ i = 0}^{ m - 1}  R (\tau A_N)^{m-i} P_N \Phi \Delta W^H_i.
\end{equation}
Consequently,
\begin{equation}\label{eq:time-error-decomposition}
\begin{split}
   \| X^N ( t_{m} ) - & \bar{ X}^N_{m} \|_{L^2 ( \Omega; V) }
     \leq
    \big\| \big( E_N ( t_m ) - R (\tau A_N)^m \big) P_N \xi \big \|_{L^2 ( \Omega; V) }
    \\ & 
    \quad +
   \Big \| \sum_{i = 0}^{m-1} \int_{t_i}^{t_{i+1} }   \big[ E_N (t_m - s ) P_N F ( X^N( s ) ) 
            -  R (\tau A_N)^{m-i} P_N F ( \bar{X}^N_{i} ) \big] \, \dd s \Big \|_{L^2 ( \Omega; V) }
    \\ & 
    \quad +
     \Big \| \sum_{i = 0}^{m-1} \int_{t_i}^{t_{i+1} }  \big[ E_N( t_m - s ) - R (\tau A_N)^{m-i} \big] P_N  \Phi \, \dd W^H(s) \Big \|_{L^2 ( \Omega; V) } 
    \\
    & \quad := \mathbb{J}_1 + \mathbb{J}_2 + \mathbb{J}_3.
\end{split}
\end{equation}
%
%
As a direct consequence of \eqref{eq:rational-exp-differ} with $\mu = \nu = 2 H + \beta - 1$,
\begin{equation} \label{eq:estimate-J1}
\mathbb{J}_1 \leq  C \tau^{ \frac{ 2 H + \beta - 1}{2} } \| \xi \|_{L^2(\Omega; V_{2 H + \beta - 1})}.
\end{equation}
%
%
In order to bound $\mathbb{J}_2$, we split it into four terms as follows:
\begin{equation} \label{eq:estimate-J2}
\begin{split}
\mathbb{J}_2 
& \leq 
 \Big \| \sum_{i = 0}^{m-1} \int_{t_i}^{t_{i+1} }   E_N (t_m - s ) P_N \big( F ( X^N( s ) ) 
            -  F ( X^N( t_i ) )  \big) \, \dd s 
 \Big \|_{L^2 ( \Omega; V) }
 \\ & \quad +
  \Big \| \sum_{i = 0}^{m-1} \int_{t_i}^{t_{i+1} }  \big( E_N (t_m - s ) 
            -  E_N (t_{m-i} ) \big) P_N F ( X^N( t_i ) )   \, \dd s 
 \Big \|_{L^2 ( \Omega; V) }
 \\ & \quad +
  \Big \| \sum_{i = 0}^{m-1} \int_{t_i}^{t_{i+1} }   \big( E_N (t_{m-i} ) 
            -  R (\tau A_N)^{m-i}  \big) P_N F ( X^N( t_i ) )  \, \dd s 
 \Big \|_{L^2 ( \Omega; V) }
 \\ & \quad + 
 \Big \| \sum_{i = 0}^{m-1} \int_{t_i}^{t_{i+1} }   R (\tau A_N)^{m-i}  P_N 
            \big( F ( X^N( t_i ) ) 
            -   F ( \bar{X}^N_{i} ) \big) \, \dd s \Big \|_{L^2 ( \Omega; V) }
\\ & :=
\mathbb{J}_{21} + \mathbb{J}_{22} + \mathbb{J}_{23} + \mathbb{J}_{24} .
\end{split}
\end{equation}
%
%
Using the stability of $E_N(t) P_N$ and \eqref{eq:X-N-optimal-time}
together with \eqref{eq:F.conditions2} yields
\begin{equation}
\mathbb{J}_{21}  
\leq C \big( 1 + \| \xi \|_{L^2 ( \Omega; V_{2 H + \beta - 1}) } \big) 
\tau^{\frac{2 H + \beta - 1}{2} }.
\end{equation} 
In view of  \eqref{eq:Semigroup-differ}, \eqref{eq:F.conditions1}, \eqref{eq:EN-P-N-commun} 
and \eqref{eq:X-N-optimal-space},  one can show that
\begin{equation} \label{eq:J22-estimate}
\begin{split}
\mathbb{J}_{22} & \leq
\sum_{i = 0}^{m-1} \int_{t_i}^{t_{i+1} }  
         \big \| \big( E_N (t_m - s ) 
            -  E_N (t_{m-i} ) \big) P_N F ( X^N( t_i ) )   \big \|_{L^2 ( \Omega; V) }  \, \dd s 
\\ & \leq
\sum_{i = 0}^{m-1} \int_{t_i}^{t_{i+1} } 
          C (t_m - s )^{ - \frac{2 H + \beta - 1}{2} } ( s - t_i )^{ \frac{2 H + \beta - 1}{2} } 
          \big \| P_N F ( X^N( t_i ) ) \big \|_{L^2 ( \Omega; V) }  \, \dd s
\\ & \leq
C \big( 1 + \| \xi \|_{L^2 ( \Omega; V ) } \big) \tau^{\frac{2 H + \beta - 1}{2} }.
\end{split}
\end{equation}
Making use of \eqref{eq:rational-exp-differ} with $\mu = 2H + \beta - 1$ and $\nu = 0$ gives
\begin{equation} \label{eq:J23-estimate}
\mathbb{J}_{23} \leq \sum_{i = 0}^{m-1} \int_{t_i}^{t_{i+1} } 
          C t_{m-i}^{ - \frac{2 H + \beta - 1}{2} } \tau^{ \frac{2 H + \beta - 1}{2} } 
          \big \| P_N F ( X^N( t_i ) ) \big \|_{L^2 ( \Omega; V) }  \, \dd s
          \leq
          C \big( 1 + \| \xi \|_{L^2 ( \Omega; V) } \big) \tau^{\frac{2 H + \beta - 1}{2} }.
\end{equation}
At last, the stability of $R (\tau A_N)P_N$ and \eqref{eq:F.conditions2} lead us to
\begin{equation} \label{eq:J24-estimate}
\mathbb{J}_{24} \leq C \tau \sum_{i = 0}^{m-1}
          \big \| X^N( t_i )  - \bar{X}^N_{i} \big \|_{L^2 ( \Omega; V) }.
\end{equation}
Putting the above estimates together results in
\begin{equation}\label{eq:J2-estimate-final}
\mathbb{J}_2 \leq  
C \big( 1 + \| \xi \|_{L^2 ( \Omega; V_{2 H + \beta - 1}) } \big)
\tau^{\frac{2 H + \beta - 1}{2} } 
+
C \tau \sum_{i = 0}^{m-1}
          \big \| X^N( t_i )  - \bar{X}^N_{i} \big \|_{L^2 ( \Omega; V) }.
\end{equation}
In the next step, we come to the estimate of $\mathbb{J}_3$. Define $\lfloor t \rfloor_\tau : = 
t_i, t \in [t_i, t_{i+1})$, $i = 0, 1,...,M-1,$ and  a continuous version of $R(\tau A_N)^j$ by
\begin{equation}
\hat{E}_N ( t ) : =  R(\tau A_N)^j, \qquad t \in [t_{j-1}, t_j), \: j = 1, 2,..., M.
\end{equation}
This together with It\^o's isometry implies that
\begin{equation}  \label{eq:J3-estimate}
\small
\begin{split} 
 | \mathbb{J}_3 |^2  
  & =
\Big \| \sum_{i = 0}^{m-1} \int_{t_i}^{t_{i+1} }  \big[ E_N( t_m - s ) - \hat{E}_N ( t_m - s )  \big] P_N  \Phi \, \dd W^H(s) 
\Big \|_{L^2 ( \Omega; V) }^2
\\ & \leq
2 \Big \| \sum_{i = 0}^{m-1} \int_{t_i}^{t_{i+1} }  \big[ E_N( t_m - s ) - E_N( t_m - \lfloor s \rfloor_\tau )  \big] P_N  \Phi \, \dd W^H(s) 
\Big \|_{L^2 ( \Omega; V) }^2
\\ & \quad +
2\Big \| \sum_{i = 0}^{m-1} \int_{t_i}^{t_{i+1} }  \big[ E_N( t_m - \lfloor s \rfloor_\tau ) - \hat{E}_N ( t_m - s )  \big] P_N  \Phi \, \dd W^H(s) 
\Big \|_{L^2 ( \Omega; V) }^2
\\ & =  
2 \int_0^{t_m} \!\! \int_0^{t_m} \Big \langle  \big[ E_N( t_m - u ) - E_N( t_m - \lfloor u \rfloor_\tau )  \big] P_N  \Phi,  
\big[ E_N( t_m - v ) - E_N( t_m - \lfloor v \rfloor_\tau )  \big] P_N  \Phi \Big \rangle_{ \mathcal{L}_2 } 
\\ & \quad \qquad \qquad \times
\phi( u - v ) \, \dd u \dd v 
\\ & \quad +  
2 \int_0^{t_m} \!\! \int_0^{t_m} \Big \langle  \big[ E_N( t_m - \lfloor u \rfloor_\tau ) - \hat{E}_N ( t_m - u )  \big] P_N  \Phi,  
\big[ E_N( t_m - \lfloor v \rfloor_\tau ) - \hat{E}_N ( t_m - v )  \big] P_N  \Phi \Big \rangle_{ \mathcal{L}_2 } 
\\ & \quad \qquad \qquad \times
\phi( u - v ) \, \dd u \dd v 
\\ & =
2 \sum_{i, j= 0}^{m-1}  
\int_{t_j}^{t_{j+1} } \!\!  \int_{t_i}^{t_{i+1} }  \Big \langle  \mathcal{S}_N (u, t_i) P_N  \Phi,  
\mathcal{S}_N (v, t_j) P_N  \Phi \Big \rangle_{ \mathcal{L}_2 } 
\phi( u - v ) \, \dd u \dd v
\\ & \quad
+ 2 \sum_{i, j= 0}^{m-1}  
\int_{t_j}^{t_{j+1} } \!\!  \int_{t_i}^{t_{i+1} }  \Big \langle  
\mathcal{T}_{N}(i) P_N  \Phi,  
\mathcal{T}_{N}(j) P_N  \Phi \Big \rangle_{ \mathcal{L}_2 } 
\phi( u - v ) \, \dd u \dd v
\\ & :=
\mathbb{J}_{31} + \mathbb{J}_{32},
%
\end{split}
\end{equation}
where for brevity we denote
\begin{equation}\label{eq:Error-operator-S-T}
\mathcal{S}_N ( u, t_i ) : = E_N( t_m - u ) - E_N ( t_m - t_i ),
\quad
\text{ and }
\quad
\mathcal{T}_{N}(i): = E_N( t_m - t_i ) - R(\tau A_N)^{m-i}.
\end{equation}
%
%
%
Before proceeding further with the estimates of the two terms in \eqref{eq:J3-estimate}, 
we should establish the following auxiliary lemmas.
\begin{lem} \label{lem:Semigroup-difference}
For $\beta \in ( 1 - 2H, 1]$,  it holds that, for any $x \in V$,
\begin{equation} \label{eq:Semigroup-difference}
\begin{split}
\sum_{i, j = 0 }^{m-1}  
\int_{t_j}^{t_{j+1} } \int_{t_i}^{t_{i+1} } 
\Big \langle  \mathcal{S}_N ( u, t_i ) P_N  x,  
 \mathcal{S}_N ( v, t_j )  P_N  x \Big \rangle  \phi( u - v ) \, \dd u \dd v
\leq 
C \tau^{2 H + \beta - 1}  \| A_N^{ \frac{ \beta - 1} {2} } P_N x \|^2.
\end{split}
\end{equation}
\end{lem}
{\it Proof of Lemma \ref{lem:Semigroup-difference}.}
Thanks to \eqref{eq:Error-operator-S-T} and Lemma \ref{lem:optimal-regularity-lemma} by setting $\delta = H$,  
one can arrive at
\begin{equation}
\small
\begin{split}
& \sum_{i, j = 0 }^{m-1}  
\int_{t_j}^{t_{j+1} } \!\! \int_{t_i}^{t_{i+1} } 
\Big \langle \mathcal{S}_N ( u, t_i ) P_N  x,  
\mathcal{S}_N ( v, t_j ) P_N  x \Big \rangle  
\phi( u - v ) \, \dd u \dd v
\\ & \quad =
\sum_{i, j = 0 }^{m-1}   \sum_{n = 1}^N
\int_{t_j}^{t_{j+1} } \!\! \int_{t_i}^{t_{i+1} } 
\big ( e^{ - \lambda_n ( t_m - u ) } -  e^{ - \lambda_n ( t_m - t_i ) } \big )
\big ( e^{ - \lambda_n ( t_m - v ) } -  e^{ - \lambda_n ( t_m - t_j ) } \big )
\\ & \qquad\qquad\qquad\qquad \times
\langle P_N x, e_n \rangle^2 \phi( u - v ) \, \dd u \dd v
\\ & \quad \leq
C \sum_{i, j = 0 }^{m-1}  \sum_{n = 1}^N
\int_{t_j}^{t_{j+1} }  \!\! \int_{t_i}^{t_{i+1} } 
e^{ - \lambda_n ( t_m - u ) } ( u - t_i)^{ H + \frac{\beta - 1}{2} } 
e^{ - \lambda_n ( t_m - v ) }  ( v - t_j)^{ H + \frac{\beta - 1}{2} }
\lambda_n^{2 H + \beta - 1} 
\\ & \qquad\qquad\qquad\qquad \times
\langle P_N x,  e_n \rangle^2 \phi( u - v ) \, \dd u \dd v
\\ &  \quad \leq
C \tau^{2 H + \beta - 1} 
\sum_{i, j = 0 }^{m-1}  \sum_{n = 1}^N
\int_{t_j}^{t_{j+1} } \!\! \int_{t_i}^{t_{i+1} } 
e^{ - \lambda_n ( t_m - u ) }  e^{ - \lambda_n ( t_m - v ) } 
\lambda_n^{2 H } \Big \langle A_N^{ \frac{ \beta - 1} {2} } P_N x, e_n \Big \rangle^2 \phi( u - v ) \, \dd u \dd v
\\ & \quad =
C \tau^{2 H + \beta - 1} 
\sum_{n = 1}^N
\int_{0}^{t_{m} } \!\! \int_{0}^{t_{m} } 
e^{ - \lambda_n ( t_m - u ) }  e^{ - \lambda_n ( t_m - v ) } 
\lambda_n^{2 H } \Big \langle A_N^{ \frac{ \beta - 1} {2} } P_N x, e_n \Big \rangle ^2 \phi( u - v ) \, \dd u \dd v
\\ & \quad =
C \tau^{2 H + \beta - 1} 
\int_{0}^{t_{m} } \!\! \int_{0}^{t_{m} } 
\Big \langle A_N^H E_N ( t_m - u ) A_N^{ \frac{ \beta - 1} {2} } P_N x, 
A_N^H  E_N ( t_m - v ) A_N^{ \frac{ \beta - 1} {2} } P_N x  \Big \rangle
\phi( u - v ) \, \dd u \dd v
\\ &  \quad \leq
C \tau^{2 H + \beta - 1}  \| A_N^{ \frac{ \beta - 1} {2} } P_N x \|^2,
\end{split}
\end{equation}
as required. Here the elementary inequality $\lambda^{ - \alpha} (1 - e^{-\lambda s} )\leq 
C s^{\alpha}, \lambda, s >0$ was also used.
$\square$
\begin{lem} \label{lem:BE-Wiener-error-estimate}
For $\beta \in ( 1 - 2H, 1]$,  it  holds that,
\begin{equation} \label{eq:BE-Wiener-error-estimate}
\begin{split}
\sum_{i, j= 0}^{m-1}  
\int_{t_j}^{t_{j+1} } \int_{t_i}^{t_{i+1} }  \Big \langle  \mathcal{T}_{N}(i) P_N  x,  
\mathcal{T}_{N}(j) P_N  x \Big \rangle  \phi( u - v ) \, \dd u \dd v
\leq 
C \tau^{2 H + \beta - 1}  \| A_N^{ \frac{ \beta - 1} {2} } P_N x \|^2.
\end{split}
\end{equation}
\end{lem}
{\it Proof of Lemma \ref{lem:BE-Wiener-error-estimate}.}
Analogously as before, we decompose the error as follows:
\begin{equation}
\small
\begin{split}
& \sum_{i, j = 0 }^{m-1}  
\int_{t_j}^{t_{j+1} } \!\! \int_{t_i}^{t_{i+1} } 
\Big \langle  \mathcal{T}_{N}(i) P_N  x,  
\mathcal{T}_{N}(j) P_N  x \Big \rangle  \phi( u - v ) \, \dd u \dd v
\\ & \quad =
\sum_{i, j = 0 }^{m-1}   \sum_{ \tau \lambda_n \leq 1 }
\int_{t_j}^{t_{j+1} } \!\!  \int_{t_i}^{t_{i+1} } 
\big ( e^{ - \lambda_n ( t_m - t_i ) } -  R(\tau \lambda_n)^{m-i}  \big )
\big ( e^{ - \lambda_n ( t_m - t_j ) } -  R(\tau \lambda_n)^{m-j}  \big )
\\ & \qquad\qquad\qquad\qquad\qquad\qquad \times
\langle P_N x, e_n \rangle^2 \phi( u - v ) \, \dd u \dd v
\\ & \qquad +
\sum_{i, j = 0 }^{m-1}   \sum_{ \tau \lambda_n > 1 }
\int_{t_j}^{t_{j+1} } \!\!  \int_{t_i}^{t_{i+1} } 
\big ( e^{ - \lambda_n ( t_m - t_i ) } -  R(\tau \lambda_n)^{m-i}  \big )
\big ( e^{ - \lambda_n ( t_m - t_j ) } -  R(\tau \lambda_n)^{m-j}  \big )
\\ & \qquad\qquad\qquad\qquad\qquad\qquad \times
\langle P_N x, e_n \rangle^2 \phi( u - v ) \, \dd u \dd v
\\ & \quad  := \mathbb{III}_1 + \mathbb{III}_2. 
\end{split}
\end{equation}
With the aid of \eqref{eq:rational-appro-error}, we start the estimate of $\mathbb{III}_1$:
\begin{equation}\label{eq:estimate-K1}
\begin{split}
| \mathbb{III}_1 | \leq & \sum_{i, j = 0 }^{m-1}   \sum_{ \tau \lambda_n \leq 1 }
\int_{t_j}^{t_{j+1} } \!\!  \int_{t_i}^{t_{i+1} } 
C ( m - i) ( m - j) ( \tau \lambda_n )^4
e^{ - c ( m - i - 1 ) \tau \lambda_n }
e^{ - c ( m - j - 1 ) \tau \lambda_n }
\\ & \qquad\qquad\qquad\qquad\quad \times
\langle P_N x, e_n \rangle^2 \phi ( u - v) \, \dd u \dd v
\\ = &
C  \tau^4
\sum_{i, j = 0 }^{m-1}   \sum_{ \tau \lambda_n \leq 1 }
\lambda_n ^ { 5 - \beta }
\int_{t_j}^{t_{j+1} } \!\!  \int_{t_i}^{t_{i+1} } 
( m - i) ( m - j) 
e^{ - c ( m - i - 1 ) \tau \lambda_n }
e^{ - c ( m - j - 1 ) \tau \lambda_n }
\\ & \qquad\qquad\qquad\qquad\qquad \times
\Big\langle A_N^{\frac{\beta - 1}{2} } P_N x, e_n \Big \rangle^2 \phi ( u - v) \, \dd u \dd v
\\ = &
C  \tau^4
\sum_{i, j = 1 }^{m}   \sum_{ \tau \lambda_n \leq 1 }
\lambda_n ^ { 5 - \beta }
\int_{ t_{m - j} }^{ t_{m - j + 1} } \! \! \int_{ t_{m - i} }^{ t_{m - i + 1} } 
i j 
e^{ - c ( i - 1 ) \tau \lambda_n }
e^{ - c ( j - 1 ) \tau \lambda_n }
\\ & \qquad\qquad\qquad\qquad\qquad \times
\Big \langle A_N^{\frac{\beta - 1}{2} } P_N x, e_n \Big \rangle^2 
\phi ( u - v ) \, \dd u \dd v
\\ = &
C  \tau^4
\sum_{i, j = 1 }^{m}   \sum_{ \tau \lambda_n \leq 1 }
\lambda_n ^ { 5 - \beta }
\int_{ t_{j -1} }^{ t_{ j } } \! \! \int_{ t_{ i - 1} }^{ t_{ i } }
i j 
e^{ - c ( i - 1 ) \tau \lambda_n }
e^{ - c ( j - 1 ) \tau \lambda_n }
\Big\langle A_N^{\frac{\beta - 1}{2} } P_N x, e_n \Big\rangle^2 
\phi ( u - v ) \, \dd u \dd v.
\end{split}
\end{equation}
Accordingly, elementary facts and Lemma \ref{lem:lambda-phi-integral} help us to derive from the above estimate that
\begin{equation}
\begin{split}
|\mathbb{III}_1 | \leq &
C \tau^2
\sum_{i, j = 1 }^{m}   \sum_{ \tau \lambda_n \leq 1 }
\lambda_n ^ { 5 - \beta } e^{ 2 c \tau \lambda_n }
\int_{ t_{j -1} }^{ t_{ j } } \!\!  \int_{ t_{ i - 1} }^{ t_{ i } }
( u + \tau ) ( v + \tau ) 
e^{ - c u \lambda_n }
e^{ - c v \lambda_n }
\\ & \qquad\qquad\qquad\qquad\qquad \times
\Big\langle A_N^{\frac{\beta - 1}{2} } P_N x, e_n \Big \rangle^2 
\phi ( u - v ) \, \dd u \dd v
\\ \leq &
C \tau^2
\sum_{ \tau \lambda_n \leq 1 }
\lambda_n ^ { 5 - \beta }
\Big \langle A_N^{\frac{\beta - 1}{2} } P_N x, e_n \Big \rangle^2 
\int_{ 0 }^{ t_{ m } } \! \! \int_{ 0 }^{ t_{ m } }
( u + \tau ) ( v + \tau ) 
e^{ - c u \lambda_n }
e^{ - c v \lambda_n }
\phi ( u - v ) \, \dd u \dd v
\\ \leq &
C \tau^{2 H + \beta - 1} 
\sum_{ \tau \lambda_n \leq 1 }
\Big \langle A_N^{\frac{\beta - 1}{2} } P_N x, e_n \Big \rangle^2. 
\end{split}
\end{equation}
Observing $ (1 + x)^{-1} \geq e^{ - x }, x > 0 $ and $ i + \tfrac{\beta - 1}{2} >  i - H > 0, i \in \N$, 
changing variables $(u - t_i )/ \tau = \bar{u}, 
(v - t_i )/ \tau = \bar{v}$, and exploiting Lemma \ref{lem:phi-estimate} one acquires
\begin{equation}\label{eq:estimate-K1}
\begin{split}
| \mathbb{III}_2 | \leq & \sum_{i, j = 0 }^{m-1}   \sum_{ \tau \lambda_n > 1 }
\int_{t_j}^{t_{j+1} } \!\!  \int_{t_i}^{t_{i+1} } 
\frac{1}{ ( 1 + \tau \lambda_n )^{ m - i } } \frac{1}{ ( 1 + \tau \lambda_n )^{ m - j } }
\langle P_N x, e_n \rangle^2 \phi ( u - v) \, \dd u \dd v
\\ = &
\tau^{ \beta - 1 } \sum_{i, j = 0 }^{m-1}   \sum_{ \tau \lambda_n > 1 }
\int_{t_j}^{t_{j+1} } \!\!  \int_{t_i}^{t_{i+1} } 
\frac{ ( \tau \lambda_n )^{1 - \beta} }{ ( 1 + \tau \lambda_n )^{ m - i } } \frac{1}{ ( 1 + \tau \lambda_n )^{ m - j } }
\big \langle A_N ^{ \frac{\beta - 1}{2} } P_N x, e_n \big \rangle^2 \phi ( u - v) \, \dd u \dd v
\\ \leq &
\tau^{ \beta - 1 } \sum_{i, j = 0 }^{m-1}   \sum_{ \tau \lambda_n > 1 }
\int_{t_j}^{t_{j+1} } \!\!  \int_{t_i}^{t_{i+1} } 
2^{- ( m - i + \frac{\beta - 1}{2} ) } 2^{- ( m - j + \frac{\beta - 1}{2} ) }
\big \langle A_N ^{ \frac{\beta - 1}{2} } P_N x, e_n \big \rangle^2 \phi ( u - v) \, \dd u \dd v
\\ = &
\tau^{ 2H + \beta - 1 } 
\sum_{ \tau \lambda_n > 1 }
\big \langle A_N ^{ \frac{\beta - 1}{2} } P_N x, e_n \big \rangle^2
\sum_{i, j = 0 }^{m-1}  
\int_0^1 \!\! \int_0^1 2^{- ( m - i + \frac{\beta - 1}{2} ) } 2^{- ( m - j + \frac{\beta - 1}{2} ) }
\phi ( \bar{u} + i - \bar{v} - j ) \, \dd \bar{u}  \dd \bar{v}
\\ = &
\tau^{ 2H + \beta - 1 } 
\sum_{ \tau \lambda_n > 1 }
\big \langle A_N ^{ \frac{\beta - 1}{2} } P_N x, e_n \big \rangle^2
\sum_{i, j = 1 }^{m}  
\int_0^1 \!\! \int_0^1 2^{- ( i + \frac{\beta - 1}{2} ) } 2^{- ( j + \frac{\beta - 1}{2} ) }
\phi ( \bar{u} + i - \bar{v} - j ) \, \dd \bar{u} \dd  \bar{v}
\\ \leq &
\tau^{ 2H + \beta - 1 } 
\sum_{ \tau \lambda_n > 1 }
\big \langle A_N ^{ \frac{\beta - 1}{2} } P_N x, e_n \big \rangle^2
\sum_{i, j = 1 }^{m}  
2^{- ( i + \frac{\beta - 1}{2} ) } 2^{- ( j + \frac{\beta - 1}{2} ) }
\max(i, j)^{2H - 1}
\\ \leq &
C \tau^{ 2H + \beta - 1 } 
\sum_{ \tau \lambda_n > 1 }
\big \langle A_N ^{ \frac{\beta - 1}{2} } P_N x, e_n \big \rangle^2.
\end{split}
\end{equation}
Gathering the above two estimates together yields the desired assertion.
$\square$

%
Now we can proceed to estimate terms in \eqref{eq:J3-estimate}.
As a direct result of Lemma \ref{lem:Semigroup-difference} and \eqref{eq:defn-L2}, 
one can infer that
\begin{equation}
\begin{split}
\mathbb{J}_{31} 
& =
2 \sum_{i, j = 0 }^{m-1}  
\int_{t_j}^{t_{j+1} } \!\! \int_{t_i}^{t_{i+1} } \Big \langle  \mathcal{S}_N ( u, t_i ) P_N  \Phi,  
\mathcal{S}_N ( v, t_j ) P_N  \Phi \Big \rangle_{ \mathcal{L}_2 } 
\phi( u - v ) \, \dd u \dd v
\\ &  \leq
C \tau^{2 H + \beta - 1}  \sum_{k \in \N} \| A_N^{ \frac{ \beta - 1} {2} } P_N \Phi e_k \|^2
\\ &  \leq 
C \tau^{2 H + \beta - 1}  \| A^{ \frac{\beta -1}{2} } \Phi \| _{\mathcal{L}_2}^2.
\end{split}
\end{equation}
Similarly, employing Lemma \ref{lem:BE-Wiener-error-estimate} enables us to derive
\begin{equation}
\small
\begin{split}
\mathbb{J}_{32}
=
2 \sum_{i, j = 0 }^{m-1} 
\int_{t_j}^{t_{j+1} } \!\! \int_{t_i}^{t_{i+1} } 
\Big \langle  \mathcal{T}_{N}(i) P_N  \Phi,  
\mathcal{T}_{N}(j) P_N  \Phi \Big \rangle_{ \mathcal{L}_2 }  
\phi( u - v ) \, \dd u \dd v
\leq
C \tau^{2 H + \beta - 1}  \| A^{ \frac{\beta -1}{2} } \Phi \| _{\mathcal{L}_2}^2.
\end{split}
\end{equation}
Plugging the above two estimates into \eqref{eq:J3-estimate} yields
\begin{equation} \label{eq:J3-final}
 | \mathbb{J}_3 |^2  \leq C \tau^{2 H + \beta - 1}  \| A^{ \frac{\beta -1}{2} } \Phi \| _{\mathcal{L}_2}^2.
\end{equation}
Combining \eqref{eq:estimate-J1}, \eqref{eq:J2-estimate-final}, \eqref{eq:J3-final} we deduce from
\eqref{eq:time-error-decomposition} that
\begin{equation}
\| X^N ( t_{m} ) - \bar{ X}^N_{m} \|_{L^2 ( \Omega; V) }
     \leq
C \big( 1 + \| \xi \|_{L^2 ( \Omega; V_{2 H + \beta - 1}) } \big) \tau^{\frac{2 H + \beta - 1}{2} } 
+
C \tau \sum_{i = 0}^{m-1}
          \big \| X^N( t_i )  - \bar{X}^N_{i} \big \|_{L^2 ( \Omega; V) }.
\end{equation} 
Applying the discrete version of the Gronwall inequality shows the desired error bound
and the proof is thus complete. $\square$ 

\section{Numerical results}
\label{sec:numerical-result}
\begin{figure}[htp]
\centering
      \includegraphics[width=3in,height=2.8in] {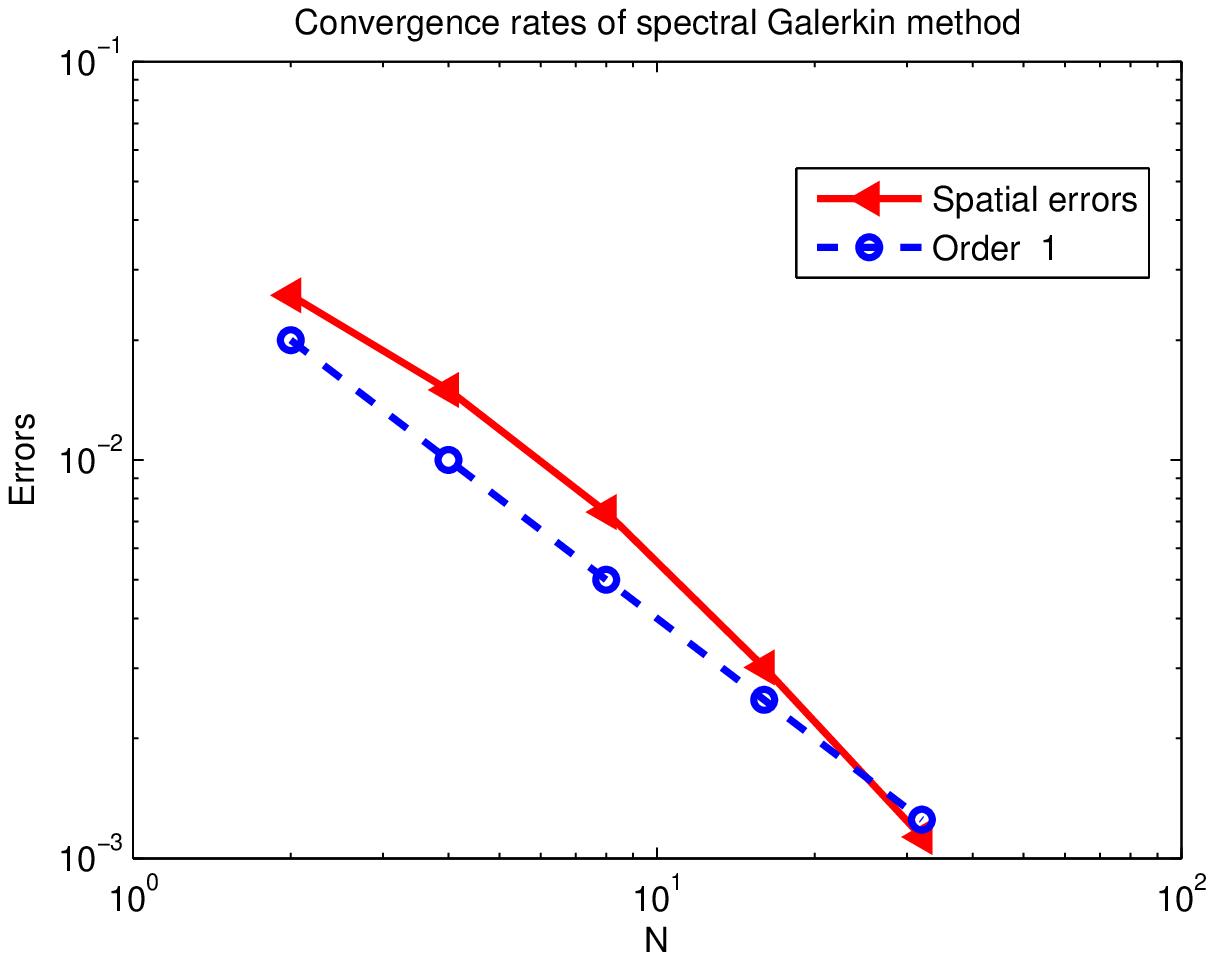}
       \includegraphics[width=3in,height=2.8in] {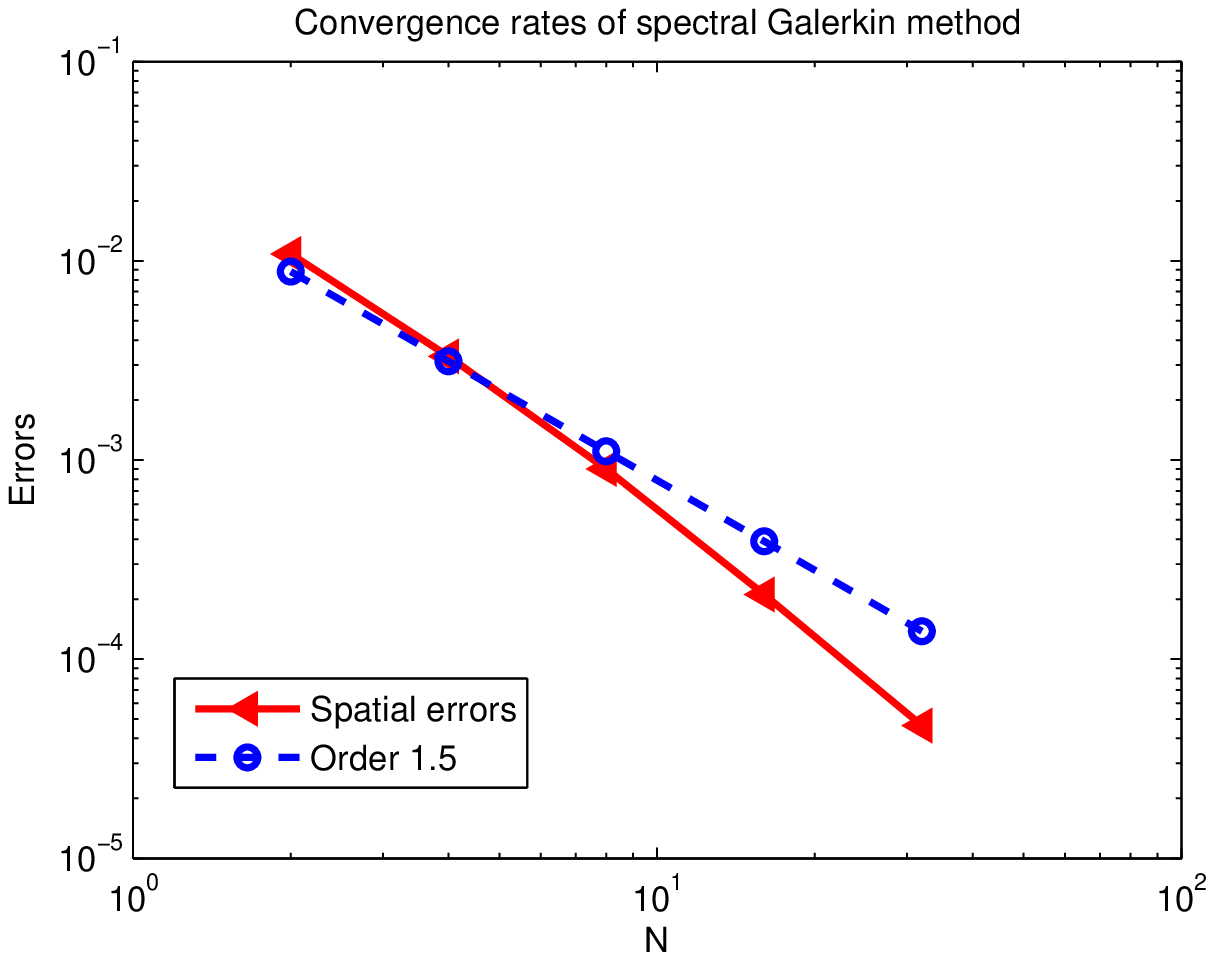}
 \caption{Convergence rates for the spatial discretization (Left: $Q=I$; right: $\text{Tr}(Q)< \infty$).}
\label{fig:spatial-error}
\end{figure}

\begin{figure}[htp]
\centering
      \includegraphics[width=3in,height=2.8in] {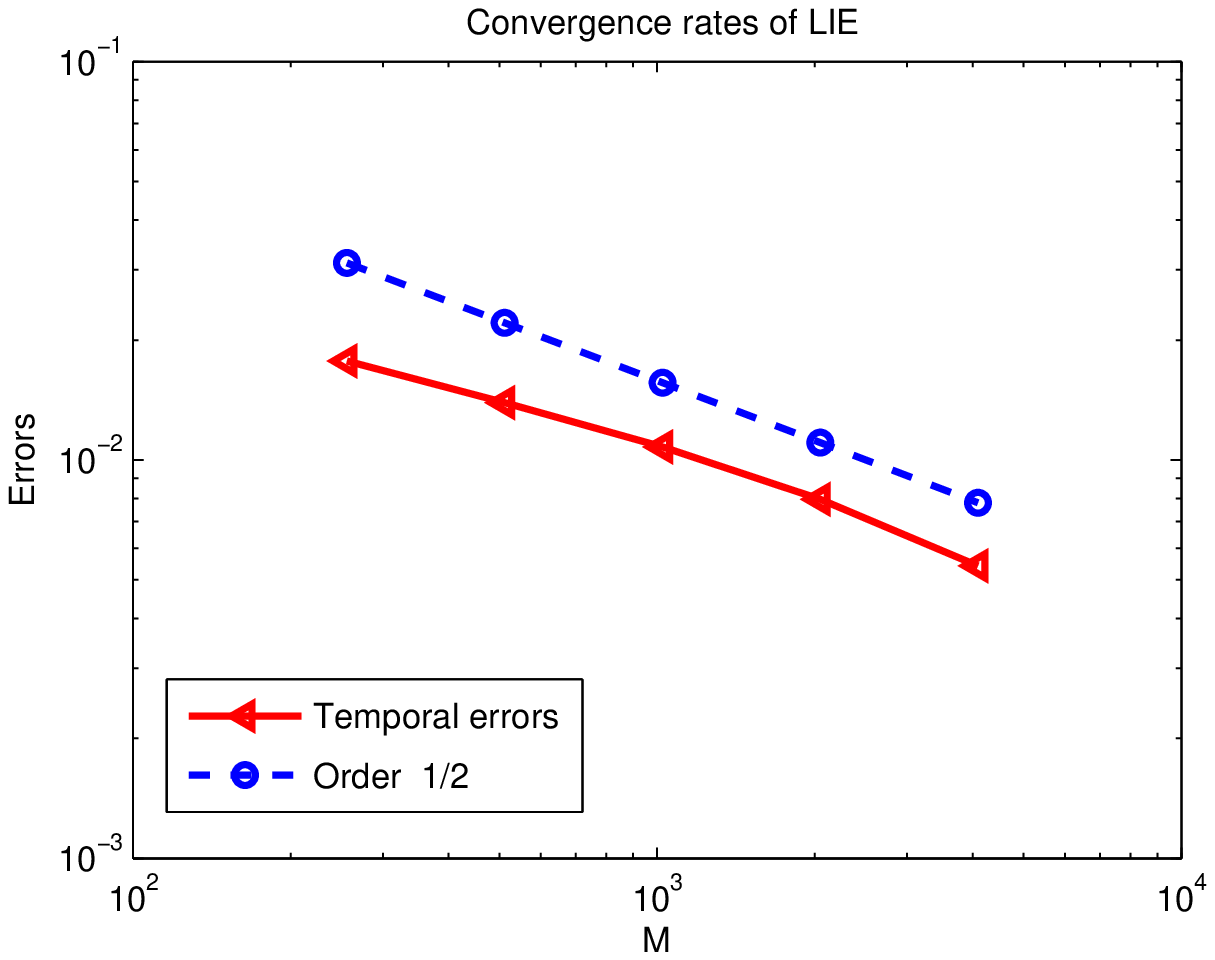}
       \includegraphics[width=3in,height=2.8in] {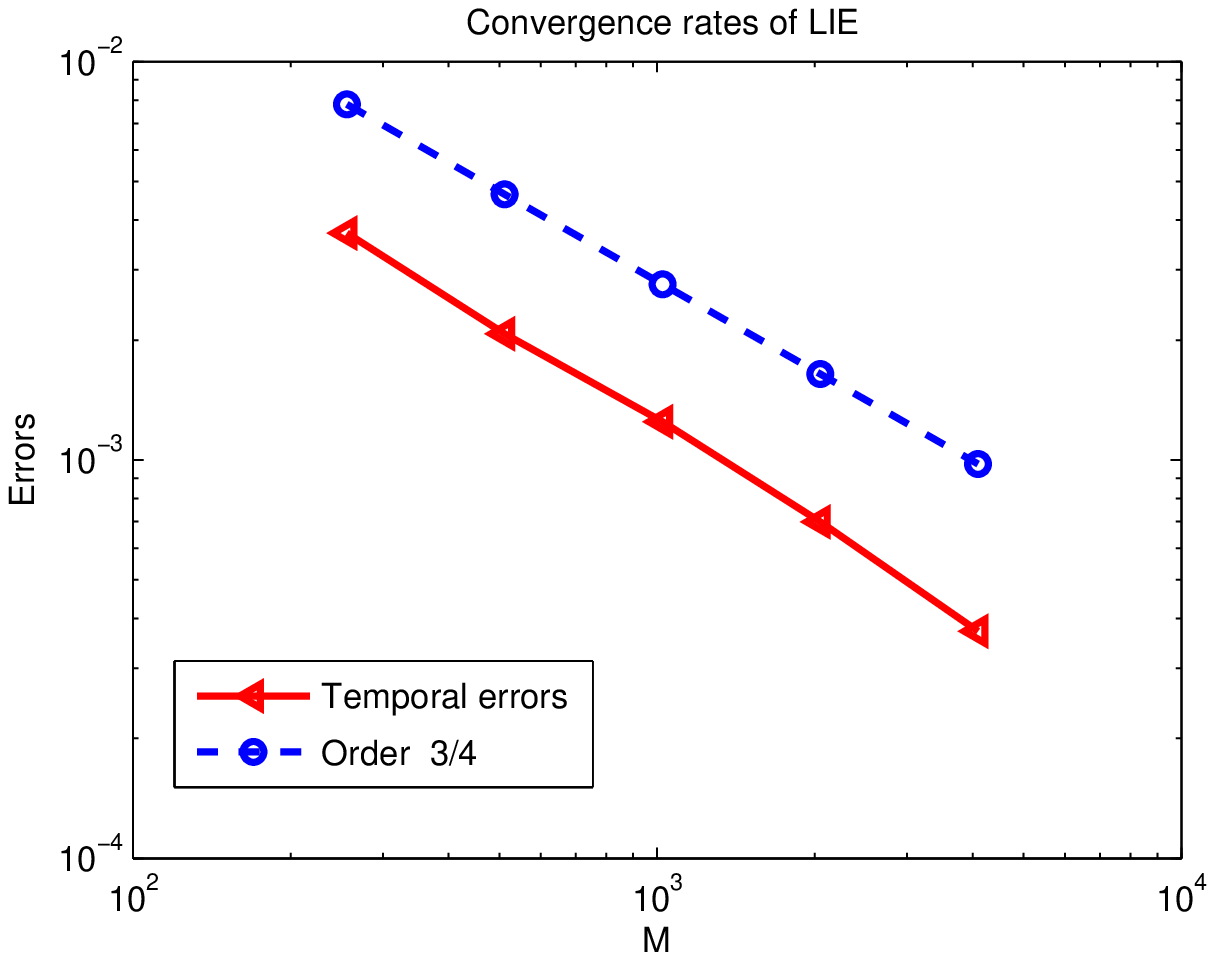}
 \caption{Convergence rates for the temporal discretization (Left: $Q=I$; right: $\text{Tr}(Q)< \infty$).}
\label{fig:temporal-error}
\end{figure}

In this section, some numerical experiments are performed to illustrate previous mean-square convergence rates.  Let us look at a test problem described by
\begin{equation}\label{eq:SHE}
\left\{
    \begin{array}{lll}
    \frac{\partial u}{\partial t} = \frac{\partial^2 u}{\partial x^2} + \sin(u) + Q^{\frac12} \dot{W}^H (t), \quad  t \in (0, 1], \:\: x \in (0,1),
    \\
     u(0, x) = \sin( \pi x), \quad\quad x \in (0,1),
     \\
     u(t, 0) = u(t,1) = 0,  \quad t \in (0, 1].
    \end{array}\right.
\end{equation}
Here $V = L^2((0, 1), \R)$ and the covariance operator  $Q \colon V  \rightarrow V$ is a bounded, linear, 
positive self-adjoint operator with a unique positive square root $Q^{\frac12}$.
In what follows we fix $H = \tfrac34$ and just consider two types of covariance operators, one being $Q = I_V$ 
and the other given by
\begin{equation} \label{eq:numer-section-Q-trace}
Q e_1 = 0, \quad Q e_i = \frac{1}{ i \log (i)^2} e_i  \quad \forall \, i \geq 2.
\end{equation}
Obviously, \eqref{eq:numer-section-Q-trace} guarantees $\text{Tr}(Q) = \| Q^{\frac12} \|_{\mathcal{L}_2} < \infty$
and thus condition \eqref{AQ_condition} is fulfilled with $\beta = 1$. When $Q = I_V$, \eqref{AQ_condition}
is then satisfied with $\beta < \frac12$. For the above two cases, we present in Fig.\ref{fig:spatial-error} and 
Fig.\ref{fig:temporal-error}  (log-log scale) mean-square approximation errors at the endpoint $T = 1$, 
caused by spatial and  temporal discretizations, respectively. 
The fBm is simulated in the sprit of  \cite{abry1996wavelet} and
the expectations are approximated by computing averages over 100 samples.

To demonstrate the convergence rates in space, we identify the ``exact'' solution by using the full 
discretization with $\tau_{\text{exact}} = 1/200$, $N_{\text{exact}} = 2^{12}$.
The spatial approximation errors $\| X(1) - X^N(1) \|_{L^2 ( \Omega; V) }$ with
$ N= 2^i, i = 1, 2, ..., 5 $ are depicted in Fig.\ref{fig:spatial-error}, where one can detect
different numerical performances for the above two different covariance operators.  
More specifically, the resulting spatial errors decrease at slopes close to $1$ and $1.5$ 
for $Q=I_V$ and $Q$ given by \eqref{eq:numer-section-Q-trace}, respectively. 
This is consistent with the previous theoretical result \eqref{eq:spatial-error}.
Likewise, we illustrate the convergence rates of the temporal approximation errors in Fig.\ref{fig:temporal-error}.
Based on the spectral Galerkin spatial discretization \eqref{eq:spectral-spde} with $N = 2^{7}$, 
we take $\tau_{\text{exact}} = 2^{-14}$ to compute the ``exact'' solution 
and five different time step-sizes $\tau = 2^{-i}, i = 8,9,...,12$ to get time discretizations. 
From Fig.\ref{fig:temporal-error},  one can observe expected convergence rates for the two distinct choices of $Q$,
which agrees with orders identified in Theorem \ref{thm:time-convergence}.

\section{Concluding remarks}
\label{sec:conclusion}
In the present paper, we attempt to provide sharp $L^2$-regularity results for 
semi-linear parabolic SPDEs driven by additive fractional noise.
In addition, mean-square convergence rates of a full discrete scheme for  
the underlying problem are studied, with optimal convergence rates in both space and time achieved.  
Much remains to be done for our future research in this area.  For example, it is still an open problem to 
recover an optimal convergence rate for the finite element spatial discretization. 
Also, sharp analysis and approximations of SPDEs driven by fBm 
with Hurst parameter $H \in (0, \tfrac12)$, as well as SPDEs driven by multiplicative fractional noise
are on the list of our future research topics.

%

\bibliographystyle{abbrv}

\bibliography{../bib/bibfile}

\end{document}